\documentclass[reqno,12pt]{amsart}

\usepackage{graphicx, amsmath, amssymb, amscd, amsthm, euscript, amsfonts,color}
\usepackage[utf8]{inputenc}
\usepackage[english]{babel}
\usepackage{verbatim}

\usepackage[pdfencoding=auto, psdextra]{hyperref}
\usepackage{bookmark}
\usepackage{psfrag}

\newtheorem{theorem}{Theorem}[section]

\newtheorem{proposition}[theorem]{Proposition}
\newtheorem{corollary}[theorem]{Corollary}
\newtheorem{claim}{Claim}[theorem]
\newtheorem{Claim}[theorem]{Claim}

\theoremstyle{remark}
\newtheorem{remark}[theorem]{Remark}

\theoremstyle{definition}
\newtheorem{definition}[theorem]{Definition}

\numberwithin{equation}{section}

\DeclareMathOperator{\Ad}{{\rm Ad}}          

\DeclareMathOperator{\Lie}{{\rm Lie}}        

\newcommand{\field}[1]{\mathbb{#1}}
\newcommand{\R}{\field{R}}
\newcommand{\N}{\field{N}}
\newcommand{\Q}{\field{Q}}
\newcommand{\Z}{\field{Z}}
\newcommand{\T}{\mathbb{T}}

\def\cF{\mathcal{F}}
\def\cP{\mathcal{P}}
\def\cH{\mathcal{H}}
\def\cL{\mathcal{L}}

\newcommand{\lieg}{\mathfrak{g}}
\newcommand{\lief}{\mathfrak{f}}
\newcommand{\liet}{\mathfrak{t}}

\providecommand{\abs}[1]{\lvert#1\rvert}
\providecommand{\Abs}[1]{\Bigl\lvert #1 \Bigr\rvert}
\providecommand{\norm}[1]{\lVert#1\rVert}

\newcommand{\inv}{^{-1}}
\newcommand{\cl}[1]{\overline{#1}}

\begin{document}

\title[Equidistribution of dilated curves]{Equidistribution of dilated curves \\ on Nilmanifolds}
\author{Bryna Kra}
\address{Northwestern University, Evanston, IL 60208 USA}
\email{kra@math.northwestern.edu}
\author{Nimish A. Shah}
\address{The Ohio State University, Columbus, Ohio 43210 USA}
\email{shah@math.osu.edu}
\author{Wenbo Sun}
\address{The Ohio State University, Columbus, Ohio 43210 USA}
\email{sun.1991@osu.edu}

\thanks{The first author was partially supported by NSF grant DMS-1500670 and the second author was partially supported by NSF grant DMS-1700394.}

\begin{abstract}
 Generalizing classic results for a family of measures in the torus, for a family $(\mu_t)_{t\geq 0}$ of measures defined on a nilmanifold $X$, we study conditions under which the family equidistributes, meaning conditions under which the measures $\mu_t$ converge  as $t\to\infty$ in the weak$^\ast$ topology to the Haar measure on $X$. We give general 
 conditions on a family of measures defined by a dilation  process, showing necessary and sufficient conditions for equidistribution as the family dilates, along with conditions such that this holds for all dilates outside some set of density zero.  Furthermore, we show that these two types of equidistribution are different.

 %
\end{abstract}    
     
\maketitle

\section{Introduction}
\subsection{Limiting distributions of measures}
A classic problem for billiards is of illumination: in a polygonal room, a light source is located at some point and the question is if there is some point not illuminated by this source. Chaika and Hubert~\cite{CH}
recently studied a related problem for circles of light, rather than points, showing 
that dilated circles around a fixed point weakly equidistribute outside a set of density zero (in their terminology, 
this phenomenon is called {\em weak illumination}).  
Motivated by their work, we prove weak and strong equidistribution results for a dilated family of measures on a nilmanifold. 

The study of equidistribution results on a nilmanifold originates in the work 
of Green~\cite{Green:homogeneous}, where he showed that a flow is either equidistributed 
or there is a nontrivial obstruction to this flow arising from horizontal character on the nilmanifold.
More precise distributional results were obtained by
Shah~\cite{Shah-polynomial}, who described the limiting behavior for a polynomial series of iterates in a flow,  
 and by Leibman~\cite{Leibman}, who proved convergence results for polynomial sequences of iterates.  
The asymptotic behavior of dilates of a measure supported on a curve in a nilmanifold was studied by 
Bj\"{o}rklund  and Fish~\cite{BF}, and among other results we answer a conjecture of theirs about the general behavior of such dilates (this result is given in Theorem~\ref{thm:analytic}).

For example, we consider a family of measures that are linear expansions on the Lie algebra $\lieg$ associated to some 
nilpotent Lie group $G$.  A sample result is necessary and sufficient conditions that this family 
equidistributes as it dilates, and whether or not these equidistribution
 results for all sufficiently large dilates 
or only for all dilates outside a set of density zero depends on the derivative of the dilation. 
The simplest case of our results is for a torus, where such strong equidistribution results for a curve 
are implicit in the literature (though we are unaware of an explicit result like this).  For a nilmanifold $X$, 
if we consider a continuous curve with the induced measure on this curve, it is easy to check that 
the dilates equidistribute weakly so long as the curve contains no linear segment of positive length (this 
corresponds to a nontrivial horizontal character picking up positive measure). More general equidistribution results 
on a nilmanifold, covering 
a broader class of dilations, require significantly more work and these characterizations are our main focus.

To give the precise formulations of our results in Sections~\ref{sec:weak-equi}, \ref{sec:equi-analytic}, and \ref{sec:equi-smooth},  we start by defining the objects that give us sufficient (and in some cases necessary) conditions for equidistribution.   

\subsection{Equidistribution and weak equidistribution}
Let $X=G/\Gamma$ be a compact nilmanifold, meaning that 
$G$ is a simply connected nilpotent Lie group and $\Gamma\subset G$ is a co-compact closed subgroup of $G$.  The group  
$G$ acts on $X=G/\Gamma$ by left translation, and there is a unique $G$-invariant Borel probability measure $\mu$ on $X$ (the Haar measure).

A family of probability measures $(\mu_{t})_{t\geq 0}$ on $X$ 
is \emph{equidistributed on $X$} if $\mu_{t}$ converges in the weak$^\ast$ topology
as $t\to\infty$  to the Haar measure $\mu$ of $X$.
Letting $C(X)$ denote the space of continuous functions on $X$, this means that the family of measures $(\mu_t)_{t\geq 0}$ is equidistributed if for all $f\in C(X)$, 
\begin{equation} \label{eq:equidistribution}
\lim_{t\to\infty}\int_{X}f\, d\mu_{t}=\int_{X}f\, d\mu. 
\end{equation}

Throughout, we let $\lambda$ denote the Lebesgue measure.  
The family of measures  $(\mu_{t})_{t\geq 0}$ is said to be \emph{weakly equidistributed on $X$} if the convergence in~\eqref{eq:equidistribution}
holds for $t$ in a set of asymptotic density $1$ in $[0,\infty)$, meaning there exists $A\subset [0,\infty)$ such that
\begin{gather*}
\lim_{T\to\infty} \frac{\lambda(A\cap [0,T])}{T}=1
\end{gather*}
and for all $f\in C(X)$, 
\begin{equation*}
\lim_{\substack{t\to\infty\\t\in A}}\int_{X}f\, d\mu_{t}=\int_{X}f\, d\mu. 
\end{equation*}
If we consider discrete time, then the family $(\mu_n)_{n\in\N}$ of measures 
is weakly equidistributed on $X$ if the convergence in~\eqref{eq:equidistribution} holds for a set of parameters $n$ 
with asymptotic density $1$, meaning along a set $A\subset\N$ satisfying $\lim_{N\to\infty} \frac{|A\cap \{1,\ldots, N\}|}{N}=1$, 
where $|\cdot|$ denotes the cardinality of the set. 

Equivalently, $(\mu_t)_{t\geq 0}$ 
 is weakly equidistributed if and only if for every $f\in C(X)$ with  $\int_{X}f\, d\mu=0$,
\begin{equation} \label{eq:weak-Illu}
\lim_{T\to\infty}\frac{1}{T}\int_{0}^{T}
\Abs{\int_{X}f\, d\mu_{t}}^{2}\, dt=0, 
\end{equation}
and  for discrete iterates $(\mu_n)_{n\in\N}$, this becomes 
\begin{equation}
\lim_{N\to\infty}\frac{1}{N} \sum_{n=0}^N \Abs{\int_X f\,d\mu_n}^2 = 0.
\label{eq:weak-Illu-discr}
\end{equation}

While it is clear that strong equidistribution implies weak equidistribution, for a family of measures in a nilmanifold these two notions are not equivalent.  Various examples illustrating the difference are given in  Section~\ref{sec:counter}. 

\subsection{Dilation associated to the Lie algebra}
Assume $G$ is a nilpotent Lie group and $\lieg$ is the associated Lie algebra. 
 Let $\Ad\colon G\to {\rm Aut}(\lieg)$ denote the adjoint representation and $\Ad(G)$ the image of $G$ under this representation.  
Note that the image $\Ad(\Gamma)$ is a Zariski dense subgroup of $\Ad(G)$ (see, for example,~\cite[Chapter 2]{R-book}). Since $\Gamma$ normalizes its connected component of the identity, which we denote by $\Gamma^0$, we have that $\Gamma^0$ is a normal Lie subgroup of $G$.  
Since $\Gamma^0$ acts trivially on $G/\Gamma$, we have that
 $G$ acts on $X=G/\Gamma$ via $G/\Gamma^0$; in other words, we can say that $G/\Gamma^0$ acts on $X=G/\Gamma$ as follows:
\[
(g\Gamma^0)\cdot x=gx \quad \text{ for all $g\in G$ and $x\in X$.}
\]

Let $\lieg$ denote the Lie algebra associated to the nilpotent Lie group $G$. Then  $\lieg/\Lie(\Gamma^0)$ is the Lie algebra associated to $G/\Gamma^0$.
We define a family of {\em dilations\/} $(\rho_t)_{t\in\R}$ such that $\rho_t\colon \lieg\to \lieg/\Lie(\Gamma^0)$ is a linear transformation, and each matrix entry of $\rho_t$ with respect to any bases of $\lieg$ and $\lieg/\Lie(\Gamma^0)$ is a polynomial in $t$. In other words, for some $m\geq 0$ and for some linear maps $B_j\colon \lieg\to \lieg/\Lie(\Gamma^0)$, $0\leq j\leq m$, we have  $\rho_t=\sum_{j=0}^m t^j B_j$ for all $t\in\R$.  

\subsection{Dynamics of measures under dilations} 
Let $\nu$ be a probability measure on $\lieg$ and let $x_{0}\in G/\Gamma$. 
For $t\geq 1$, let $\mu_{\nu,x_{0},\rho_t}$ denote the measure on $X$ defined as follows: for any $f\in C(X)$, 
\begin{equation}
    \label{eq:define-measure}
\int_{X} f\, d\mu_{\nu,x_{0},\rho_t} = \int_{\lieg} f(\exp\circ \rho_t(y)\cdot x_{0})\,d\nu(y),
\end{equation}
where $\exp\colon \lieg/\Lie(\Gamma^0)\to G/\Gamma^0$ is the exponential map from Lie algebra to Lie group. 

Motivated by~\cite{BF,Randol}, we consider the following special case of $\nu$ in Sections~\ref{sec:equi-analytic} and~\ref{sec:equi-smooth}: given a measurable map $\phi\colon (0,1)\to \lieg$, we can take $\nu$ to be the pushforward of the Lebesgue measure $\lambda$ restricted to $(0,1)$ under $\phi$. 
Then for $t\geq 0$, taking $\mu_{\phi,x_{0},\rho_t}$ to be the measure on $X$ defined by
\begin{equation}\label{eq:define-measure2}
\int_{X} f\,d\mu_{\phi,x_{0},\rho_t}=\int_{0}^{1}f(\exp\circ\rho_t\circ\phi(u)\cdot x_{0})\,d\lambda(u)
\end{equation}
for all $f\in C(X)$, we have that  $\mu_{\phi,x_{0},\rho_{t}}=\mu_{\nu,x_{0},\rho_t}$.

\subsection{Abelianization of a nilsystem} 
Let $[G,G]$ denote the Lie subgroup of $G$ corresponding to the commutator algebra $[\lieg,\lieg]$.
If $\Gamma$ is  a co-compact, closed subgroup of $G$, 
then $[G,G]\Gamma$ is a closed subgroup of $G$, and $\bar X:=G/[G,G]\Gamma$ can be identified with a compact torus $\R^{m}/\Z^{m}$ for some $m\geq 1$ (see for example~\cite[Corollary~1 of Theorem~2.3]{R-book}).    
We call $\bar X$ the {\em abelianization} of $X$ and use $\liet\cong \R^m$ to denote the associated Lie algebra of $\bar X$ and  $q\colon X=G/\Gamma \to \bar X = G/[G,G]\Gamma$ to denote the natural quotient map. 

Let $dq\colon \lieg/\Lie(\Gamma^0) \to \liet$ denote the differential of the map $q$ at the identity coset. Then 
\begin{equation} \label{eq:dq}
q(\exp(y)\cdot x_0)=\exp\circ dq(y)+q(x_0)
\end{equation}
for all $y\in\lieg/\Lie(\Gamma^0)$  and all $x_0\in X$. 

Furthermore, $dq\circ\rho_t\colon \lieg\to \liet$ can be expressed as
\begin{equation} \label{eq:torus}
dq\circ\rho_t=\sum_{i=0}^{d_1} t^i A_i,
\end{equation}
where $d_1\in\N$ and each $A_i\colon \lieg\to \liet$ is a linear map.

When $\Gamma$ is discrete, and so $\Gamma^0=\{e\}$, a natural choice for a family of dilations is to consider $(\rho_t)_{t\in\R}$, where $\rho_tv=tv$ for all $v\in \lieg$. In this case, $d_1 = 1$, $A_1=dq$,  and $A_0=0$.

Let $\bar X^\ast$ denote the space of (continuous) unitary characters on the torus $\bar X$. Then elements of $\bar X^\ast$ are in one-to-one correspondence with unitary characters on $G$ such that the kernel 
of the map contains $\Gamma$. For any $\chi\in \bar X^\ast$, let $d\chi\colon\liet\to \R$ denote the differential of $\chi$, 
meaning that  for $y\in\liet$, 
		\begin{equation} \label{eq:dchi}
		    \chi(y+\Z^m)=e^{2\pi i d\chi(y)}.  
		\end{equation}

\subsection{Weak equidistribution of dilated measures}
\label{sec:weak-equi}
Our first result provides necessary and sufficient conditions for weak equidistribution of 
the family $(\mu_t:=\mu_{\nu,x_{0},\rho_t})_{t\geq 0}$ given  in~\eqref{eq:define-measure}: 
\begin{theorem}\label{thm:willumination}
Let $\nu$ be a probability measure on $\lieg$, $x_0\in X$, and $(\mu_t:=\mu_{\nu,x_{0},\rho_t})_{t\geq 0}$ be the family of measures defined in~\eqref{eq:define-measure}.  If 
for all $\chi\in \bar X^\ast\setminus\{1\}$ and $(z_1,\ldots,z_{d_1})\in\R^{d_1}$ we have 
\begin{equation} \label{eq:hyp2}
		    \nu(\{v\in\lieg\colon d\chi(A_iv)=z_i \text{ for all  } 1\leq i\leq d_1\})=0, 
\end{equation}
then the families of measures $(\mu_t)_{t\geq 0}$ and $(\mu_n)_{n\in\N}$ are weakly equidistributed on $X$. 
		
Furthermore, if $A_{0}=0$, then condition~\eqref{eq:hyp2} is necessary for weak equidistribution of $(\mu_n)_{n\in\N}$, and also for $(\mu_t)_{t\geq 0}$. 
\end{theorem}

In Remark~\ref{rem:willumination-necesity}, we provide an example exhibiting why we assume $A_0=0$ in proving necessity in  Theorem~\ref{thm:willumination}. 

\subsection{Equidistribution of dilated analytic curves}
\label{sec:equi-analytic}
We maintain the same notation in this section.   
Carrying out the constructions  of the family of measures in~\eqref{eq:define-measure2} 
with an analytic map $\phi\colon (0,1)\to\lieg$, we 
resolve a conjecture stated in Bj\"{o}rklund  and Fish~\cite{BF} (see the discussion following Theorem~7 in their paper): 

\begin{theorem} \label{thm:analytic} 
Assume $\phi\colon (0,1)\to\lieg$ is analytic, let $x_0\in X$, and let $(\mu_t=\mu_{\phi,x_0,\rho_{t}})_{t\in\R}$ be as defined in~\eqref{eq:define-measure2}. 
If for every $\chi\in\bar X^\ast\setminus\{1\}$, there exists $1\leq i\leq d_1$ such that the map 
\begin{equation}
 \label{eq:analytic-cond}
u\mapsto d\chi(A_i\phi(u)) 
\end{equation}
is not constant, then $(\mu_t)_{t\geq 0}$ is equidistributed on $X$. 

Furthermore, if $A_0=0$, then the conditions on the maps given in~\eqref{eq:analytic-cond} are 
also necessary for the equidistribution of $(\mu_t)_{t\geq 0}$ on $X$.
\end{theorem}

\subsection{Equidistribution of dilated differentiable curves} 
\label{sec:equi-smooth}
Our main result provides a condition on tangents to the curve for (strong) equidistribution, when $\phi$ is a sufficiently differentiable curve. Still maintaining the assumptions and notation stated at the beginning of this section, we have: 
\begin{theorem} \label{thm:smooth}
Assume that $X = G/\Gamma$ is a compact nilmanifold. There exists a natural number $D$, which can be expressed in terms of the degrees of polynomials in $t$ defining $\rho_t$, such that the following holds: suppose that $\phi^{(D)}(u)$ exists for (Lebesgue) almost all $u\in(0,1)$ and  that 
for every $\chi\in \bar X^\ast\setminus \{1\}$ and (Lebesgue) almost every $u\in (0,1)$,
\begin{equation} \label{eq:tangent}
 A_i\phi^{(1)}(u)\notin \ker d\chi
\end{equation}
for some $1\leq i\leq d_1$.
Then the family of measures $(\mu_{\phi,x_{0},\rho_t})_{t\geq 0}$ is equidistributed on $G/\Gamma$.
\end{theorem}
Our proof is based on a stronger equidistribution result for dilates of curves that shrink (the precise statement is given in Theorem~\ref{thm:local:illumination}). 

The following special case is of interest: 
\begin{corollary} \label{cor:torus} Let $X=\R^n/\Z^n$ and $\phi\colon (0,1)\to \R^n$ be an almost everywhere differentiable map. Suppose that for each $\bf v\in \Z^n\setminus\{0\}$, 
$\phi^{(1)}(u)$ is not orthogonal to ${\bf v}$ for (Lebesgue) almost every $u\in (0,1)$.
Then for any $f\in C(\R^n/\Z^n)$ and $x_0\in\R^n/\Z^n$,
\[
\lim_{t\to\infty}\int_{0}^{1} f(t\phi(s)+x_0)\,d\lambda(s)=\int_{\R^m/\Z^m} f\,d\mu,
\]
where $\lambda$ denotes the Lebesgue measure on $[0,1]$ and $\mu$ the Haar measure on $\R^n/\Z^n$. 
\end{corollary}

For analytic $\phi$, the conclusion of Corollary~\ref{cor:torus} was obtained in Randol~\cite{Randol}.

In Section~\ref{sec:counter} we  provide examples where condition in~\eqref{eq:hyp2} holds and so $(\mu_t)_{t\geq 0}$ is weakly equidistributed on $X$, but $(\mu_t)_{t\geq 0}$ is not equidistributed on $X$.  

\section{Weak equidistribution for measures}

This section is devoted to proving Theorem~\ref{thm:willumination}.  We give a detailed proof for the discrete case, 
and then indicate the modifications needed for the continuous setting.  

The classic equidistribution theorem in the discrete setting is due to Weyl: 
\begin{theorem}[Weyl~\cite{Weyl}] \label{thm:Weyl}
Let $p(t)=a_{d}t^{d}+\dots+a_{1}t+a_{0}$ be a real polynomial. Then 
\[
\lim_{N\to\infty} \frac{1}{N}\sum_{n=1}^N e^{2\pi i p(n)} =0
\]
if and only if at least one of $a_1,\ldots,a_d$ is not rational.
\end{theorem}

If $G$ is a nilpotent Lie group, $a_{1},\dots,a_{m}\in G$, and
 $p_{1},\dots,p_{m}\colon\N\to\N$ are polynomials taking integer values on the integers, then a sequence $(g(n))_{n\in\N}$
 in $G$
 of the form $a^{p_{1}(n)}_{1}a^{p_{2}(n)}_{2}\dots a^{p_{m}(n)}_{m}$ is a \emph{polynomial}. 
 If $\lieg$ denotes the Lie algebra of $G$, we
  say that $\zeta\colon\N\to \lieg$ is a \emph{polynomial map} if $\exp(\zeta(n))=g(n)$
for some polynomial $(g(n))_{n\in\N}$.  
Generalizing Weyl's equidistribution result, Leibman~\cite{Leibman} showed (we write his result in terms of the Lie algebra): 
\begin{theorem}[Leibman~\cite{Leibman}]
	\label{thm:Leibman}
	Let $\zeta\colon\N\to \lieg$ be a polynomial map such that for every nontrivial character $\chi$ on $\bar X$, 
	\begin{equation} \label{eq:rationaldifference}
	d\chi\circ dq\circ \zeta(n)\not\equiv d\chi\circ dq\circ \zeta(0)\mod \Q \text{ for some $n\in\N$}.
\end{equation}
Then for any $f\in C(X)$ and any $x\in X$,
	\begin{equation*}
	\lim_{N\to\infty} \frac{1}{N}\sum_{n=1}^{N} f(\exp\circ \zeta(n)\cdot x)=\int_X f\,d\mu.
	\end{equation*}
\end{theorem}

We make use of these results in the proof of the weak equidistribution result: 
\begin{proof}[Proof of Theorem~\ref{thm:willumination} (in the discrete setting)]
First we prove the necessity of condition~\eqref{eq:hyp2} for weak equidistribution. 
Assume that $A_0=0$ and suppose that~\eqref{eq:hyp2} fails to hold, and we want to prove that the families $(\mu_t)_{t\geq 0}$ and $(\mu_n)_{n\in\N}$ are both not weakly equidistributed. 
Thus there exist a nontrivial character $\chi$ on $\bar X$ and $z_{1},\dots,z_{d_{1}}\in\R$ such that 
\begin{equation} \label{eq:uchi}
\nu(\{v\in \lieg\colon d\chi(A_{j} v)=z_{j} \text{ for all } 1\leq j\leq d_{1}\})
\end{equation}
is positive.  
Since $\chi$ is a nontrivial character on the torus $\bar X$ and its pushforward
$q_\ast\mu$ is the translation invariant (Haar) probability measure on the compact torus $\bar X$, we have

\begin{equation} \label{eq:zero}
   \int_{X} \chi\circ q\, d\mu=\int_{\bar X} \chi \,d(q_\ast\mu)=0.
\end{equation}

Therefore in view of~\eqref{eq:weak-Illu-discr}, in order prove that $(\mu_n)_{n\in\N}$ is not weakly equidistributed, we need to show that
\[
\limsup_{N\to\infty}\frac{1}{N}\sum_{n=1}^{N}\Abs{\int_{X} \chi\circ q\, d\mu_n}^2>0.
\]
Therefore it suffices to prove that

\begin{equation} \label{eq:Mfactorial}
    \lim_{M\to\infty}\lim_{N\to\infty}\frac{1}{N}\sum_{n=1}^{N}\Abs{\int_{X} \chi\circ q\, d\mu_{M! n}}^2>0,
\end{equation}
because for any sequence $\{a_n\}$ of non-negative reals, and 
$M\in\N$,
\[
\frac{1}{M!N}\sum_{n=1}^{M!N}a_n \geq \frac{1}{M!} \frac{1}{N}\sum_{n=1}^{N} a_{M! n}.
\]

For all $t\in\R$,
\begin{align*}
 \int_{X} \chi\circ q\, d\mu_t  & = \int_{\lieg} \chi(q(\exp\circ \rho_{t}(y)\cdot x_0))\,d\nu(y) \\
 &=\chi(q(x_0))\int_{\lieg} \chi(dq\circ\rho_{t}(y)+\Z^m)\,d\nu(y)\\
 &=\chi(q(x_0))\int_{\lieg} e^{2\pi i d\chi(dq\circ\rho_{t}(y))}\,d\nu(y)\\
 &=\chi(q(x_0))\int_{\lieg} 
 e^{2\pi i \sum_{j=0}^{d_{1}}t^{j}\cdot d\chi(A_{j}y)}\,d\nu(y), 
 \end{align*}
where we have used equations~\eqref{eq:dq}, \eqref{eq:dchi}, and~\eqref{eq:torus}.  
Writing $y_{j}:=d\chi(A_{j}y)$  and using the assumption that $A_0=0$, this last quantity is the same as
$$
 \chi(q(x_0))\int_{\lieg} 
 e^{2\pi i \sum_{j=1}^{d_{1}}t^{j}\cdot y_{j}}\,d\nu(y).  
$$
Combining this expression for $ \int_{X} \chi\circ q\, d\mu_t$ with the fact that $\abs{\chi(q(x_0))}=1$, we have that
 \begin{align*}
 &\lim_{M\to\infty} \lim_{N\to\infty}\frac{1}{N}\sum_{n=1}^{N} 
    \Abs{\int_{X} \chi\circ q\, d\mu_{M!n}}^2  \\
 &=\lim_{M\to\infty} \lim_{N\to\infty}\frac{1}{N}\sum_{n=1}^{N} 
    \int_{\lieg\times\lieg}
e^{2\pi i (nM!)^j\sum_{j=1}^{d_{1}}(y_{j}-y'_{j})}\,d(\nu\times\nu)(y,y')
\\
 &=
    \int_{\lieg\times\lieg}\left[\lim_{M\to\infty}\lim_{N\to\infty}\frac{1}{N}\sum_{n=1}^{N} e^{2\pi i \sum_{j=1}^{d_{1}} n^j((M!)^{j}(y_{j}-y'_{j}))}\right]d(\nu\times\nu)(y,y').
    \end{align*}
We claim that this integral over $\lieg\times\lieg$ is 
the same as integrating the indicator function of the set
\[
W = \{(y,y')\in\lieg\times\lieg\colon y_{j}\equiv y'_{j} \mod \Q \text{ for all } 1\leq j\leq d_{1}\}.
\] 
To see this, fix $y,y'\in\lieg$. If $y_{j}\neq y'_{j} \mod \Q \text{ for some } 1\leq j\leq d_{1}$, by Weyl Equidistribution (Theorem~\ref{thm:Weyl}) the limit equals to 0. If $y_{j}= y'_{j}  \mod \Q \text{ for all } 1\leq j\leq d_{1}$, then $(M!)^{j}y_{j}= (M!)^{j}y'_{j}  \mod \Z \text{ for all } 1\leq j\leq d_{1}$ and all $M$ sufficiently large, and so the limit as $N\to\infty$ equals to $1$.

Thus 
\begin{align*}
 &\quad\lim_{N\to\infty}\frac{1}{N}\sum_{n=1}^{N} 
    \Abs{\int_{X} \chi\circ q\, d\mu_n}^2  
    =\int_{\lieg\times\lieg}\textbf{1}_{W}\,d(\nu\times\nu)(y,y') 
     \\
    &=\int_{\lieg} \nu(\{y\in \lieg\colon d\chi(A_{j} y)\equiv y'_j \mod \Q \text{ for all } 1\leq j\leq d_{1}\})\, d\nu(y')
    \\
    &=\sum_{w_{1},\dots,w_{d_{1}}\in\mathbb{\R}}
    \nu(\{y\in \lieg\colon d\chi(A_{j} y)\equiv w_{j}\mod\Q \text{ for all } 1\leq j\leq d_1\})\times\\
    &\quad \quad\quad \nu(\{y'\in \lieg\colon d\chi(A_{j} y')=w_{j} \text{ for all } 1\leq j\leq d_1\})\\&
    \geq\nu(\{v\in \lieg\colon d\chi(A_{j} v)=z_{j} \text{ for all } 1\leq j\leq d_{1}\})^2.
\end{align*}
By~\eqref{eq:uchi}, this is strictly positive and so ~\eqref{eq:Mfactorial} holds, which completes the proof that $(\mu_n)_{n\in\N}$ is not weakly equidistributed on $X$.

\subsubsection*{The converse implication.} Assume that condition~\eqref{eq:hyp2} of Theorem~\ref{thm:willumination} holds. Our goal is to 
show that~\eqref{eq:weak-Illu-discr} holds for the family 
$(\mu_{n}:=\mu_{\nu,x_{0},\rho_n})_{n\in\N}$ and any given $f\in C(X)$ with $\int_{X}f\, d\mu=0$.  We have:  
\begin{align}   
&\quad\lim_{N\to\infty}\frac{1}{N}
\sum_{n=1}^{N} \Abs{\int_{X}f\, d\mu_{n}}^{2}
\nonumber \\
&=\lim_{N\to\infty}\frac{1}{N}\sum_{n=1}^{N}
\left[\int_{X\times X} 
f\otimes\overline{f} \,d(\mu_{n}\times \mu_{n})\right] 
\nonumber \\
&=\lim_{N\to\infty}\frac{1}{N}\sum_{n=1}^{N}
\left[\int_{\lieg\times \lieg} f(\exp\circ \rho_{n}(y)\cdot x_{0}) \overline{f(\exp\circ \rho_{n}(y')\cdot x_{0})} \,d(\nu\times \nu)(y,y')\right]
\nonumber\\
&=\int_{\lieg\times\lieg}
\left[\lim_{N\to\infty}\frac{1}{N}\sum_{n=1}^{N} f\otimes \bar f(\exp\circ \zeta_{y,y'}(n)\cdot(x_{0},x_{0})) \right]\,d(\nu\times\nu)(y,y'),
\label{eq:equat1}
\end{align}
where $\zeta_{y,y'}\colon\N\to \lieg/\Lie(\Gamma^0)\times\lieg/\Lie(\Gamma^0)$ is the map given by 
\begin{equation} \label{eq:zeta}
\zeta_{y,y'}(t)=(\rho_t (y),\rho_t (y')) \quad\text{ for all } t\in\N.
\end{equation}
From the definition of $(\rho_t)_{t\in\R}$, it follows that $(\exp\circ\zeta_{y,y'}(n))_{n\in\N}$ is a polynomial on $(G/\Gamma^0)\times (G/\Gamma^0)$.   

For every $(\chi,\chi')\in \bar X^\ast \times \bar X^\ast\setminus \{(1,1)\}$, we have $d(\chi\otimes \bar\chi')=d\chi\oplus d\bar\chi'$, and by~\eqref{eq:torus},
\begin{align*}
    &(d\chi\oplus d\bar\chi')\circ\zeta_{y,y'}(n)\\
    &=d\chi(dq\circ\rho_n(y))-d\chi'(dq\circ\rho_n(y'))\\ 
    &=\sum_{j=1}^{d_1} [d\chi(A_jy)-d\chi'(A_jy')]n^j + d\chi(A_0y)-d\chi'(A_0y'). 
\end{align*}

Therefore in view of condition~\eqref{eq:rationaldifference} of Theorem~\ref{thm:Leibman}, 
$$
(d\chi\oplus d\bar\chi')\circ\zeta_{y,y'}(n)\equiv d\chi(A_0y)-d\chi'(A_0y') \mod \Q 
$$
for all $n \in \Z$
is equivalent to the condition that 
$$
d\chi(A_jy)\equiv d\chi'(A_jy') \mod\Q
$$ for all $1\leq j \leq d_1$.  

In order to apply Theorem~\ref{thm:Leibman} to $\zeta_{y,y'}$ and $X\times X$,  define 
$B$ to be the set of $(y,y')\in \lieg\times\lieg$ such that 
 for all $(\chi,\chi')\in \bar X^\ast \times \bar X^\ast\setminus \{(1,1)\}$, 
there exists  $j\in\{1, \ldots, d_1\}$ such that  $d\chi(A_jy)\not\equiv d\chi(A_jy')\mod\Q$.  
Define $C$ to be the set 
 of $y\in \lieg$ such that  for all $\chi\in \bar X^\ast\setminus\{1\}$, 
 there exists $j\in\{1, \ldots, d_1\}$ such that $d\chi(A_jy)\not\in\Q$ 
and define 
$$
E_y=\{y'\in \lieg \colon (y,y')\not\in B\}.
$$

Since $\Q$ is countable, by the hypothesis given in~\eqref{eq:hyp2}, we have that $\nu(C)=1$. Also 
\begin{equation} \label{eq:CBA1}
    B\supset\bigcup_{y\in C} \{y\}\times (\lieg\setminus E_y).
\end{equation}

For any $y\in C$, $\chi\in \bar X^\ast$, $\chi'\in\bar X^\ast\setminus \{1\}$, define $E_{y,\chi,\chi'}$ to be the set
\[
\{y'\in \lieg\colon d\chi'(A_jy')\equiv d\chi(A_jy)\mod\Q \text{ for all } 1\leq j\leq d_1\}.
\]
Then since $\Q$ is countable, by the hypothesis~\eqref{eq:hyp2}, we have that 
\begin{equation} \label{eq:Eychi}
    \nu(E_{y,\chi,\chi'})=0.
\end{equation}

Suppose that $y\in C$, $y'\in E_y$, and $(\chi,\chi')\in \bar X^\ast\times \bar X^\ast\setminus\{(1,1)\}$ are such that 
\begin{equation}\label{eq:2.9}
d\chi(A_jy)\equiv d\chi'(A_jy')\mod\Q \text{ for all $1\leq j\leq d_1$}.
\end{equation}
By the definition of the set $C$, we can choose $j$ with $1\leq j\leq d_1$ such that $d\chi(A_jy)\notin \Q$. Then by~\eqref{eq:2.9}, $\chi'(A_jy')\not\in \Q$. Therefore $\chi'\neq 1$. 
By the definition of the set $E_y$, 
\begin{equation} 
E_{y}=\bigcup_{\chi'\in \bar X^\ast\setminus\{1\}} \bigcup_{\chi\in\bar X^\ast} E_{y,\chi,\chi'}.
\end{equation}
Since $\bar X^\ast$ is countable, by~\eqref{eq:Eychi}, we have $\nu(E_y)=0$. Hence by~\eqref{eq:CBA1}, 
$(\nu\times\nu)(B)=1$. 

By Theorem~\ref{thm:Leibman} applied to $X\times X$, for any $f\in C(X)$ with $\int_X f\,d\mu=0$, for all $ (y,y')\in B$, we have
\[
\lim_{N\to\infty} \frac{1}{N} \sum_{n=1}^N f\otimes \bar f(\exp\circ \zeta_{y,y'}(n)\cdot(x_{0},x_{0}))=0.
\]

Now since $(\nu\times\nu)((\lieg\times \lieg) \setminus B)=0$, it follows from~\eqref{eq:equat1} that
\begin{align*}
     &\quad\lim_{N\to\infty} \frac{1}{N} \sum_{n=1}^N \Abs{\int_\lieg f(\exp\circ \rho_{n}(y)\cdot x_{0})\,d\nu(y)}^2 \\
  &=\int_{\lieg\times\lieg} \Bigl[\lim_{N\to\infty} \frac{1}{N} \sum_{n=1}^N f\otimes \bar f(\exp\circ\zeta_{y,y'}(n)\cdot(x_{0},x_{0}))\Bigr] d(\nu\times\nu)(y,y')=0.
\end{align*}
Thus  $(\mu_n)_{n\in\N}$ is weakly equidistributed. 
\end{proof}

\subsubsection*{The case of continuous parameter} The proof in this case is similar. Instead of using the discrete version of Weyl's Theorem to show necessity, we use: 
 \begin{theorem}[Weyl~\cite{Weyl}]\label{thm:Weyl-cont}
Let $p(t)=a_{d}t^{d}+\dots+a_{1}t+a_{0}$ be a polynomial with real coefficients. Then 
$$\lim_{T\to\infty}\frac{1}{T}\int_{0}^{T}e^{2\pi ip(t)}\,dt=0$$
if and only if at at least one of $a_{1},\dots,a_{d}$ is nonzero. 
\end{theorem}

We also replace the use of Leibman's Theorem by the following result of Shah (which generalizes Theorem~\ref{thm:Weyl-cont}): 
\begin{theorem}[Shah~\cite{Shah-polynomial}]
	\label{thm:Shah}
	Let $\zeta\colon\R\to \lieg$ be a polynomial map such that for every nontrivial character $\chi$ on $\bar X$, 
	\begin{equation} \label{eq:nonconstant}
	d\chi\circ dq\circ \zeta(t)\neq d\chi\circ dq\circ \zeta(0)\quad 
	\text{ for some $t\in\R$}.
	\end{equation}
	Then for any $f\in C(X)$ and any $x\in X$,
	\begin{equation*}
	\lim_{T\to\infty} \frac{1}{T} \int_0^T f(\exp\circ \zeta(t)\cdot x)\,dt=\int_X f\,d\mu.
	\end{equation*}
\end{theorem}
Up to obvious changes in notation, the proof is then the same as in the discrete case; in fact it simplifies, as one can replace equivalence $\mod \Q$ by equality.  We omit the details, 

\begin{remark} \label{rem:willumination-necesity}
We give an example illustrating that without the assumption that $A_0=0$ 
in Theorem~\ref{thm:willumination}, 
the necessity condition for weak equidistribution becomes quite complicated. 
Let $G=\R^2$, $\Gamma=\Z^2$, $X=\R^2/\Z^2=\R/\Z\times \R/\Z$, $\rho_t=tA_1+A_0$ $t\in\R$, where $A_1(x_1,x_2)=(0,x_1)$, and $A_0(x_1,x_2)=x_1$ for all $(x_1,x_2)\in\R^2$, and let $\nu=\nu_1\times \delta_{0}$, where $\nu_1$ is the standard Lebesgue measure on $\R$ restricted to $[0,1)$ and $\delta_0$ denotes the unit mass at $0$ on $\R$. Let $\chi$ be a character on $G/\Gamma$ such that $d\chi(x_1,x_2)=x_1$ for all $(x_1,x_2)\in\R^2$. Then 
\[
\nu(\{v\in\lieg\colon d\chi(A_1v)=0\})=1.
\]
Therefore the condition~\eqref{eq:hyp2} fails to hold in this case. On the other hand, 
we will verify that for any $f\in C(X)$ and $(\mathbf{x_0},\mathbf{y_0})\in X$, we have 
\begin{equation} \label{eq:m}
\lim_{t\to\infty} \int_{X} f\,d\mu_t=\lim_{t\to\infty}\int_{0}^1 f(x+\mathbf{x_0},tx+\mathbf{y_0})\,dx =\int_{X} f\,dm, 
\end{equation}
where $m$ denotes the normalized Lebesgue measure on $X$; in other words, the family $\{\mu_t\}_{t\geq 0}$ is equidistributed on $X$, and hence weakly-equidistributed on $X$. 

To prove~\eqref{eq:m}, let $\lambda$ be a weak$^\ast$ limit of $\mu_{t_i}$ for some sequence $t_i\to\infty$. We are left with checking that $\lambda=m$.  
To see this, let $y\in\R$ and define 
$$
f_y(\mathbf{x_1},\mathbf{x_2})=f(\mathbf{x_1},y+\mathbf{x_2}) \quad \text{ for all }(\mathbf{x_1},\mathbf{x_2})\in X.
$$
Then
\begin{align}
    \int_X f_y\,d\lambda &= \lim_{i\to\infty}\int_X f_y \,d\mu_{t_i}=\lim_{i\to\infty} \int_0^1 f(x+\mathbf{x_0},y+t_ix+\mathbf{y_0})\,dx \nonumber\\
    &= \lim_{i\to\infty} \int_0^1 f(x+\mathbf{x_0},t_i(t_i^{-1}y+x)+\mathbf{y_0})\,dx\nonumber\\
    &=\lim_{i\to\infty} \int_{t_i^{-1}y}^{t_i^{-1}y+1} f(-t_i^{-1}y+x+\mathbf{x_0},t_ix+\mathbf{y_0})\,dx\nonumber\\
    &=\lim_{i\to\infty} \left[\int_0^1 f(x+\mathbf{x_0},t_ix+\mathbf{y_0})\,dx + O(2t_i^{-1}y\norm{f}_{\infty})+O(\varepsilon_{t_i})\right]\nonumber\\
    &=\int_X f\, d\lambda,
    \label{eq:f_y}
\end{align}
where $\varepsilon_{t_i}=\sup_{x\in \R}\abs{f(-t_i^{-1}y+x+\mathbf{x_0},t_ix+\mathbf{y_0})-f(x+\mathbf{x_0},t_ix+\mathbf{y_0})}\to 0$ as $i\to\infty$, because of uniform continuity of $f$.

We define $\bar f$ for $(\mathbf{x_1},\mathbf{x_2})\in X$ by  setting 
$$\bar f(\mathbf{x_1},\mathbf{x_2})  =\int_0^1 f_y(\mathbf{x_1},\mathbf{x_2})\, dy  =\int_0^1 f(\mathbf{x_1},y+\mathbf{x_2})\, dy.
$$
Then $\bar f$ is constant with respect to the second coordinate, and by \eqref{eq:f_y},
\begin{align*}
   \quad\int_X f\,d\lambda & =\int_0^1 \Bigl[\int_X f_y\, d\lambda\Bigr]\,dy
    = \int_X \bar f\,d\lambda 
     =\lim_{i\to\infty} \int_X \bar f\, d\mu_{t_i}\\
    & = \lim_{i\to\infty} \int_0^1\bar f(x+\mathbf{x_0},t_ix+\mathbf{y_0})\,dx \\ 
    &  = \lim_{i\to\infty} \int_0^1 \bar f(x+\mathbf{x_0},\mathbf{y_0})\,dx
     =\int_X f\,dm. 
\end{align*}
Thus as $t_i\to\infty$, every limit measure of the family $\{\mu_t\}_{t\geq 0}$ equals $m$, 
and so~\eqref{eq:m} holds. 
\end{remark}

\section{Strong equidistribution}
\subsection{A stronger equidistribution result for dilations of curves}

The key ingredient for proving Theorem~\ref{thm:smooth} is a stronger equidistribution of dilation of shrinking curves as stated below. We recall the setting and notation. 

We continue to use  $\lambda$ to denote the Lebesgue measure on $(0,1)$ or $\mathbb{R}$, depending on the context. 
We consider the dynamics on $X=G/\Gamma$, where $G$ is a simply connected nilpotent Lie group and $\Gamma$ is a closed subgroup of $G$. Let $\lieg$ denote the Lie algebra of $G$.

\newcommand{\gbar}{\bar\lieg}
We write $\gbar=\lieg/\Lie(\Gamma^0)$. For $t\in\R$, let $\rho_t\colon \lieg \to \gbar$ denote the linear map such that $\rho_t$ is a polynomial in $t$ with coefficients which are linear maps from $\lieg$ to $\gbar$.

Let $\liet$ denote the Lie algebra associated to the compact torus $\bar X=G/[G,G]\Gamma$. Let $dq\colon \gbar \mapsto \liet= \gbar/[\gbar,\gbar]$ denote the natural quotient map. Then $dq\circ\rho_t\colon \lieg\to \liet$ can be expressed as
\[
dq\circ\rho_t=\sum_{i=0}^{d_1} t^iA_i,
\]
where each $A_i\colon \lieg\to \liet$ is a linear map. 

Consider the lower central series
$\gbar^{(k)} = [\gbar, \gbar^{(k-1)}]$ for $k\geq 1$. Set  $\gbar^{(0)}=\gbar$ 
and let $\kappa\geq 1$ be an integer such that $\gbar^{(\kappa)}=0$.

For $1\leq k\leq \kappa$, 
let $Q_k\colon\gbar\to \gbar/\gbar^{(k)}$ denote the natural quotient map and 
choose an integer $D_k\geq 0$ such that the map $Q_k\circ\rho_t\colon \lieg\to \gbar/\gbar^{(k)}$ is a polynomial of degree at most $D_k$ in $t$ and $D_1\leq D_2\leq \cdots \leq D_{\kappa}$.  Set 
	\begin{equation} \label{eq:d}
	D=\max\left\{\sum_{m=1}^n D_{k_m}\colon  \sum_{m=1}^n k_m\leq \kappa,\, 1\leq k_m\leq \kappa, \, 1\leq n\leq \kappa\right\}.
	\end{equation}
Clearly $D\geq \kappa$. This definition of $D$ is crucially used in Claim~\ref{claim:psi}.

\begin{theorem} \label{thm:local:illumination} 
Suppose that $\phi\colon (0,1)\to\lieg$ is a function such that 
$\phi^{(D)}(u)$ exists for (Lebesgue) almost all $u\in(0,1)$. Assume further  that for every nontrivial character $\chi$ on $\bar X$ and (Lebesgue) almost every $u\in (0,1)$,
\begin{equation} \label{eq:tangent2}
A_i\phi^{(1)}(u)\notin \ker d\chi \text{ for some $1\leq i\leq d_1$.}
\end{equation}
Then there exists $W\subset (0,1)$ with 
$\lambda(W)=1$ such that for all $s_{0}\in W$ the following holds: Let $(\ell_{t})_{t\geq 1}$ be a family such that $\ell_{t}>0$, and $\ell_{t}\to \infty$ and $\ell_{t}t^{-1}\to 0$ as $t\to\infty$. Then
	\begin{equation} \label{eq:nut}
	\lim_{t\to\infty} 
	\frac{1}{\ell_t t\inv} \int_{s_0}^{s_0+\ell_t t^{-1}} f(\exp\circ \rho_t\circ \phi(\xi)\cdot x_{0})\,d\xi=\int_X f\,d\mu.
	\end{equation}
\end{theorem}

We defer the proof of this result to the next two sections. Before this, we explain its role in completing the proof Theorem~\ref{thm:smooth} and obtain some of its immediate consequences. 

\begin{proof}[Proof of Theorem~\ref{thm:smooth}]
	Let $f\in C(G/\Gamma)$ be such that $\int_{G/\Gamma} f \,d\mu=0$, 
	and  without loss of generality we can assume that $\norm{f}_\infty\leq 1$.
	Given $\varepsilon>0$, choose $W\subset (0,1)$ with $\lambda(W)=1$ as in 
	Theorem~\ref{thm:local:illumination}. The conclusion of Theorem~\ref{thm:local:illumination} can be interpreted as follows: for every $s\in W$, there exists 
	$\ell_s\geq 1$ such that for any $\ell\geq \ell_s$, there exists $t_{s,\ell}> \ell(1-s)^{-1}$ such that for all $t\geq t_{s,\ell}$,
\begin{align} \label{eq:equid}
	\Abs{\int_{s}^{s+\ell t\inv} f(\exp\circ \rho_t\circ \phi(\xi)\cdot x_{0})\,d\xi}\leq \ell t\inv \varepsilon.
\end{align}
	
	For $L>0$, let  
	\[
	W(L)=\{s\in W\colon \ell_s\leq L\}.
	\]
	By choice of $\ell_s$, there exists $L>0$ such that $\lambda(W\setminus W(L))\leq \varepsilon$.
	
	For $t\geq 1$, let 
	\[
	W(L,t)=\{s\in W(L)\colon  t_{s,L}\leq t\}.
	\]
	By choice of $t_{s,L}$, there exists $T_{L}\geq 1$ such that $\lambda(W(L)\setminus W(L,T_L))<\varepsilon$.
	
	Pick a compact set $W_2\subset W(L,T_{L})$ such that $\lambda(W(L,T_{L})\backslash W_2)\leq \varepsilon$ and let $t\geq T_L$. We cover $W_{2}$ by disjoint intervals $[s_{i},s_{i}+L t\inv]$ for $0\leq i\leq n$, where $s_0=\min W_2$, and if $s_i\in W_2$ is chosen for some $i\geq 0$, then choose
	\begin{equation*}
	s_{i+1}= \min\bigl(W_2\setminus (0,s_i+L t^{-1})\bigr)\in W_2, 
\text{ if $s_i+Lt^{-1}\leq \max (W_2)$},
\end{equation*}
otherwise put $n=i$ and stop the induction. Since $s_n\in W(L,T_L)$, and $t\geq T_L$, we have $s_n+Lt^{-1}<1$. Therefore the intervals $[s_i,s_i+Lt^{-1})$ for $0\leq i\leq n$ are disjoint and contained in $(0,1)$. Hence
\[
(n+1)(Lt^{-1})\leq 1.
\]
	
	Let $W'_{2}=\bigcup_{i=0}^{n}[s_{i},s_{i}+Lt^{-1})$.
	Then $W_{2}\subseteq W'_{2}$ and by~\eqref{eq:equid}, 
	\begin{align*}
	\Abs{\int _{W'_{2}} f(\exp\circ \rho_{t}\circ \phi(s)\cdot x_{0})\,ds} & \leq \sum_{i=0}^{n} \Abs{ \int_{s_{i}}^{s_{i}+L t\inv} f(\exp\circ \rho_{t}\circ \phi(s)\cdot x_{0})\,ds} \\ 
	& \leq (n+1)(L t\inv) \varepsilon\leq \varepsilon.  
	\end{align*}
	Since $\lambda(W_2')\geq \lambda(W_2)\geq 1-3\varepsilon$, for all $t\geq T_L$, 
	\begin{equation*}
	\Abs{\int_{0}^{1} f(\exp\circ \rho_t\circ \phi(s)\cdot x_{0})\,ds} \leq \varepsilon +3\varepsilon\norm{f}_\infty\leq 4\varepsilon.
	\qedhere
	\end{equation*}
\end{proof}

The following special case of Theorem~\ref{thm:smooth} is of interest: 

\begin{corollary} \label{cor:illumination}
	Let $\kappa\geq 1$ be such that $\gbar^{(\kappa)}=0$. Suppose that that for $1\leq k\leq\kappa$, 
	\begin{equation} \label{eq:deg}
	\text{$Q_k\circ\rho_t\colon \lieg\to \gbar/\gbar^{(k)}$ is a polynomial of degree at most $k$ in $t$}. 
	\end{equation}
	Let $\phi\colon(0,1)\to \lieg$ be such that $\phi^{(\kappa)}(u)$ exists for almost all $u\in (0,1)$. Suppose that 
	the condition on $\chi\in\bar{X}^*\setminus\{1\}$ given in~\eqref{eq:tangent} holds. 
	Then the family of measures $(\mu_{\phi,x_0,\rho_t})_{t\geq 0}$ is equidistributed on $X$ as $t\to\infty$.  
\end{corollary}

\begin{remark} \label{rem:D=kappa}
\begin{enumerate}
\item The family $\{\rho_t\}$ of automorphisms, given by each $\rho_t$ being the multiplication by $t$ on $\lieg$, is an example satisfying~\eqref{eq:deg} and generalizing Corollary~\ref{cor:torus} to nilpotent groups. 

\item Other natural examples of families of dilations $(\rho_t)_{t\in\R}$ satisfying the condition given in~\eqref{eq:deg} 
come from Stratified  Lie  algebras  and  Carnot  groups, see~\cite{BD-dilate2, BD-dilate, pansu}. 

\item In particular, if $\{\rho_t\}$ is a family of Lie algebra automorphisms of $\lieg$, then $\{\rho_t\}$ satisfies~\eqref{eq:deg}.
\end{enumerate}
\end{remark}

\begin{proof}[Proof of Corollary~\ref{cor:illumination}]
	Given the assumption in~\eqref{eq:deg}, the bound $D_k$ on the degree of the quotient map $Q_k\circ\rho_t$ 
	satisfies 
	$D_k=k$ for all $1\leq k\leq\kappa$. Therefore in~\eqref{eq:d}, we have 
	$$
	\sum_{m=1}^n D_{k_m}=\sum_{m=1}^n k_m \leq \kappa,
	$$
	and hence $D\leq \kappa$. Thus the result is a special case of Theorem~\ref{thm:smooth}.
\end{proof}

We also complete the proof of Theorem~\ref{thm:analytic}: 
\begin{proof}[Proof of Theorem~\ref{thm:analytic}]
	Since $\phi$ is an analytic curve, the condition on the $\chi$ given in~\eqref{eq:analytic-cond} implies the condition 
	given in~\eqref{eq:tangent}. Therefore the equidistribution follows from Theorem~\ref{thm:smooth}. 
	
	For the converse, suppose that the condition on the map given in~\eqref{eq:analytic-cond} fails to hold and $A_0=0$. Then there exists a nontrivial character $\chi$ on $\bar X$ such that $d\chi(A_i\phi^{(1)}(u))=0$ for all $u\in (0,1)$ and all $1\leq i\leq d_1$. But then by analyticity of $\phi$, for every $1\leq i\leq d_1$, there exists $z_i\in\R$, such that $d\chi(A_i\phi(u))=z_i$ for all $u\in(0,1)$. Hence, since $A_0=0$, 
	\begin{align}
	d\chi(dq\circ \rho_t\circ \phi(u))
	=\sum_{i=1}^{d_1} t^i d\chi(A_i\phi(u))= \sum_{i=1}^{d_1} t^i z_i. \label{eq:chiq}
	\end{align}
	Let $f=\chi\circ q\in C(X)$. Then
	\[
	\int_X f \, d\mu = \int_{\bar X} \chi\, d(q_\ast\mu)=0,
	\]
	because $\chi$ is a nontrivial character on the torus $\bar X$ and $q_\ast\mu$ is the invariant probability measure on $\bar X$.
	
	For any $x_0\in X$, by~\eqref{eq:dchi},
	\begin{equation} \label{eq:f}
	f(\exp\circ \rho_t \circ \phi(u)\cdot x_0)=e^{2\pi i d\chi(dq\circ \rho_t\circ \phi(u))}\cdot\chi(q(x_0)).
	\end{equation}	
If  $\lim_{t\to\infty}\int f\,d\mu_t$ exists, then by~\eqref{eq:chiq} and~\eqref{eq:f} $\lim_{t\to\infty} e^{2\pi i \sum_{i=1}^{d_1} t^iz_i} $ exists. We conclude that $z_i=0$ for all $i$.  
It follows that the limit is equal to $\chi(q(x_0))\neq 0$. Therefore the sequence $\mu_t$ does not converge to $\mu$ as $t\to\infty$ with respect to the weak$^\ast$ topology.
\end{proof}

\subsection{Equidistribution of expanding translates of shrinking curves}
\label{sec:Property-W}

First we describe a quantitative condition that we want for $s_0\in(0,1)$ to hold, so that~\eqref{eq:nut} holds. 

Assume that $\phi\colon(0,1)\to \lieg$ is a function satisfying the hypotheses of Theorem~\ref{thm:local:illumination}. That is, $\phi^{(D)}(u)$ exists for almost all $u\in (0,1)$, where $D$ is defined as in~\eqref{eq:d}, and for every $\chi\in\bar X^\ast\setminus\{1\}$, for almost all $u\in (0,1)$, there exists $1\leq i\leq d_1$, such that $A_i\phi^{(1)}(u)\not\in \ker d\chi$. 

Define a map $\psi_t\colon \{(u,\xi)\colon u+\xi,u\in(0,1)\}\to \lieg/\Lie(\Gamma^0)$ by
\begin{equation} \label{eq:psi}
\psi_t(u,\xi)
 =\log[\exp(\rho_t\circ\phi(u+\xi))\cdot \exp(-\rho_t\circ\phi(u))],
\end{equation}
where $\exp\colon \lieg/\Lie(\Gamma^0)\to G/\Gamma^0$ is the exponential map, and $\log$ is its inverse.

\newcommand{\cW}{\mathcal{W}}

\begin{definition}[Property $\cW(\Gamma^0)$]
We say that $s_0\in (0,1)$ has property $\cW(\Gamma^0)$ if the following holds: there exists $\alpha\geq 1$ such that for any sequences $t_n\to\infty$ and $\ell_n\to\infty$ with $\ell_{n} t_{n}^{-1}\to 0$ and for any $\varepsilon>0$, there exists a compact set $I_\varepsilon\subset [0,1]$ with $\lambda([0,1]\setminus I_{\varepsilon})<\varepsilon$ such that the following holds: for any $s\in I_\varepsilon$ and $\zeta\in \R$, for each $n\in\N$, if we put $u_n=s_0+s\ell_{n} t_{n}^{-1}$  and $\xi_n=\zeta t_n^{-\alpha}$, then
\begin{equation} \label{eq:main1}
\lim_{n\to\infty} \psi_{t_n}(u_n,\xi_n)=\eta(\zeta),
\end{equation}
where $\eta\colon \R\to\lieg/\Lie(\Gamma^0)$ is a non-constant polynomial map. 

If $\Lambda$ is any closed subgroup of $G$ containing $\Gamma$, then let $dp\colon \gbar=\lieg/\Lie(\Gamma^0)\to \lieg/\Lie(\Lambda^0)$ be the natural quotient map. Then we replace the dilations $\rho_t$ with $dp\circ \rho_t\colon \lieg \to \lieg/\Lie(\Lambda^0)$. With these modifications, property $\cW(\Lambda^0)$ is also defined. 
\end{definition}

\begin{proposition} \label{prop:local:equi}
Suppose that $s_0\in(0,1)$ has property $\cW(\Lambda^0)$ for every closed subgroup $\Lambda$ of $G$ containing $\Gamma$. For all $t\geq 1$, let $\ell_{t}>0$ be such that as $t\to\infty$, $\ell_{t}\to \infty$ and $\ell_{t}t^{-1}\to 0$. For all  $t\geq 1$ such that $s_0+\ell_tt\inv <1$, let $\nu_{t}$ be the probability measure on $X$ such that all $f\in C(X)$,
\begin{align*}
    \int_X f \,d\nu_{t} &= \frac{1}{\ell_t t\inv} \int_{s_0}^{s_0+\ell_t t^{-1}} f(\exp(\rho_t\circ \phi(\xi))\cdot x_{0})\,d\xi\\&=\int_0^1 f(\exp\circ\rho_t\circ\phi(s_0+s\ell_t t^{-1})\cdot x_0)\,ds.
\end{align*}
Then $\nu_{t}\to \mu$ with respect to the weak$^\ast$ topology as $t\to\infty$. 
\end{proposition}

\begin{proof}
The result is trivial if $\dim(G/\Gamma^0)=0$. We intend to prove the result by induction on $\dim(G/\Gamma^0)$. In particular, we can assume that the result is valid for all closed subgroups $\Lambda$ of $G$ containing $\Gamma$ such that $\dim(\Lambda^0)>\dim(\Gamma^0)$. In other words, if $p \colon  G/\Gamma \to G/\Lambda$ is the natural quotient map, and $p_\ast$ denotes the corresponding pushforward of measures then as $t\to\infty$, $p_\ast(\nu_t)$ converges to the $G$-invariant probability measure on $G/\Lambda$ with respect to the weak$^\ast$ topology. 

Let $t_{n}\to\infty$ be any given sequence. It suffices to show that
after passing to a subsequence of $\{t_n\}$, we have that $\nu_{t_n}\to \mu$ in the weak$^\ast$ topology as $n\to\infty$, where $\mu$ denotes the
$G$-invariant probability measure on $G/\Gamma$.

Since $G/\Gamma$ is compact, after passing to a subsequence, 
we may assume that $\nu_{t_{n}}\to \nu$ in the space of probability measures on $G/\Gamma$ with respect to the weak$^\ast$ topology as $n\to\infty$.

We characterize $\nu$ in stages given in the next two claims, with the result following once we can show that $\nu = \mu$. 
\begin{claim} \label{claim:invariant} 
The measure $\nu$ is invariant under a connected subgroup $U$ of $G$ properly containing $\Gamma^0$.
\end{claim}

\begin{proof}
Let $f\in C(G/\Gamma)$. Since $s_0$ has property $\cW(\Gamma^0)$, let $\alpha\geq 1$ be such that~\eqref{eq:main1} holds. For the remainder of the proof, we make 
use of some shorthand to simplify formulas. For all $n\in\N$, we write $t = t_{n}$, for all $s\in [0,1]$, we write $u=s_{0}+s\ell_{t}t^{-1}$,  and 
for all $\zeta\in\R$, we write $\xi=\zeta t^{-\alpha}$.   
Then, since $\nu_{t_n}\to \nu$ as $n\to\infty$, by~\eqref{eq:nut} we have 
\begin{align}
\int_{G/\Gamma}   f(x) \,d\nu(x)  =  \lim_{n\to\infty}
\int_{0}^{1} f(\exp\circ\rho_{t_n}\circ\phi(u)\cdot x_{0}) \,ds.    
\label{eq:piexprho}
\end{align}
Now $\xi \ell_{t}^{-1}t=\zeta \ell_t^{-1} t^{-(\alpha-1)}\to 0$, because $\alpha\geq 1$ and $\ell_t\to\infty$ as $n\to\infty$. Therefore~\eqref{eq:piexprho} can be rewritten as 
\begin{align*}
\lim_{n\to\infty}& \int_{I_{\varepsilon}+\xi \ell_{t}^{-1}t}
f(\exp\circ\rho_t\circ\phi(u)\cdot x_{0}) \,ds + O(\norm{f}_{\infty}\varepsilon)  
\\
=&\lim_{n\to\infty} \int_{I_{\varepsilon}+\xi \ell_{t}^{-1}t}
f(\exp\circ\rho_t\circ\phi(s_0+s\ell_t t^{-1})\cdot x_{0}) \,ds + O(\norm{f}_{\infty}\varepsilon)  \\
=&\lim_{n\to\infty} \int_{I_{\varepsilon}}
f(\exp\circ\rho_t\circ\phi(s_0+(s+\xi \ell_t^{-1} t)\ell_t t^{-1})\cdot x_{0}) \,ds + O(\norm{f}_{\infty}\varepsilon)  \\
=&\lim_{n\to\infty}\int_{I_{\varepsilon}}
f(\exp\circ\rho_t\circ\phi(u+\xi)\cdot x_{0})\,ds + O(\norm{f}_{\infty}\varepsilon).
\end{align*}
Applying~\eqref{eq:psi}, this becomes
\begin{align*}
\lim_{n\to\infty}  &
\int_{I_{\varepsilon}} f(\exp(\psi_{t}(u,\xi)) \exp(\rho_t\circ\phi(u))\cdot x_{0})\,ds + O(\norm{f}_{\infty}\varepsilon) \\
=& \lim_{n\to\infty}\int_{I_{\varepsilon}} f(\exp(\eta(\zeta))\exp(\rho_t\circ\phi(u))\cdot x_{0})\,ds  + O(\norm{f}_{\infty}\varepsilon)\\
=& \lim_{n\to\infty} \int_{G/\Gamma} f(\exp(\eta(\zeta))\cdot x) \,d\nu_{t_{n}}(x) + O(2\norm{f}_{\infty}\varepsilon), 
\end{align*}
where we have used~\eqref{eq:main1} in the second to last step. 
Taking the limit, in view of~\eqref{eq:piexprho}, 
$$
\Abs{\int_{G/\Gamma}   f(x) \,d\nu(x) - \int_{G/\Gamma} f(\exp(\eta(\zeta))x) \,d\nu(x)}\leq  O(2\norm{f}_{\infty}\varepsilon).
$$
Thus $\nu$ is invariant under the action of $U_1$, where $U_1$ is the closure of the subgroup of $G/\Gamma^0$ generated by $\{\exp(\eta(\zeta))\colon \zeta\in \R\}$.  
But $\zeta\mapsto \eta(\zeta)$ is a non-constant polynomial map and so  $U_1$ is a nontrivial connected subgroup of $G/\Gamma^0$. Let $U=\pi\inv(U_1)$. Then $\nu$ is $U$-invariant, $U$ is connected, and  properly contains $\Gamma^0$. 
\end{proof}

\begin{claim} \label{claim:2} The measure  $\nu$ is invariant under a closed connected normal subgroup $F$ of $G$ containing $U$ such that $F\Gamma/\Gamma$ is compact. 
\end{claim}  

\begin{proof}
Let $\cH$ be the collection of all closed connected subgroups $H$ of $G$ containing $\Gamma^0$ with $H\Gamma/\Gamma$ is compact. Then $\cH$ is countable. By Lesigne~\cite[Theorem 2]{Lesigne}, or more generally Ratner's Theorem~\cite{Ratner:solvable}, any $U$-ergodic invariant probability measure on $G/\Gamma$ is of the form $g\mu_{H}$ for some $H\in \cH$, where $\mu_{H}$ denotes the $H$-invariant probability measure on $H\Gamma/\Gamma$, and $g\in G$ is such that $Ug\Gamma\subset gH\Gamma$, or equivalently, $g\in X(U,H)$, where 
\begin{equation}
    X(U,H)=\{g\in G\colon  g\inv Ug\subset H\}.
\end{equation} 
Note that $X(U,H)$ is an algebraic subvariety of $G$ with respect to any algebraic group structure on $G$. 

Let $F\in \cH$ be a subgroup of smallest possible dimension such that 
\begin{equation} \label{eq:F}
    \nu(X(U,F)\Gamma/\Gamma)>0.
\end{equation} 
Note that
\begin{equation} \label{eq:FZ}
    X(U,H)Z=X(U,H),
\end{equation} where $Z\subset G$ is the inverse image 
of the center of $G/\Gamma^0$. Now  $\Lambda:=Z\Gamma$ is closed and $G/\Lambda$ is compact, see 
\cite[Chapter II]{R-book}. Consider the natural quotient map $p \colon  G/\Gamma \to G/\Lambda$. Let $p_\ast$ denote corresponding the pushforward map of measures. Then $p_{\ast}(\nu_{t_n})\to p_{\ast}(\nu)$ as $t_n\to\infty$. 

If $G=\Gamma^0$, the lemma is trivial. So we assume that $G/\Gamma^0$ is a nontrivial 
connected nilpotent group. Therefore its center is of strictly positive dimension. Therefore 
$\dim(G/Z)<\dim(G/\Gamma^0)$. So by our induction hypothesis stated on the beginning of the proof of Proposition~\ref{prop:local:equi}, we conclude that $p_{\ast}(\nu)$ is the $G$-invariant probability measure on $G/Z\Gamma$. 

By~\eqref{eq:F} and~\eqref{eq:FZ},
$$
 p_{\ast}(\nu)(X(U,F)Z\Gamma/Z\Gamma)=\nu(X(U,F)\Gamma/\Gamma)>0.
$$
Since $X(U,F)Z=X(U,F)$ and $p_\ast(\nu)$ is the $G$-invariant measure on $G/Z\Gamma$, the Haar measure of $X(U,F)$ is strictly positive. As we noted before, $X(U,F)$ is an algebraic subvariety of $G$. Therefore $X(U,F)=G$; that is, $Ug\subset gF$ for all $g\in G$. 

Next we show that that $F$ is normal in $G$, and $\nu$ is $F$-invariant.  
By minimality of $F$,  
$$  
\nu(S(U,F)\Gamma/\Gamma)=0, 
$$
where
$$
  S(U,F):=\bigcup\{X(U,H)\colon H\in\cH,\,\dim H<\dim F\}.
$$
Therefore 
\[
\nu(X(U,F)\Gamma/\Gamma\setminus S(U,F)\Gamma/\Gamma)=1.
\]
As a consequence, almost every $U$-ergodic component of $\nu$ is of the form $g\mu_F$ for some $g\in G$. 

Now $G\setminus S(U,F)\neq \emptyset$. 
Let $g\in G\setminus S(U,F)$. We claim that $\cl{Ug\Gamma}=gF\Gamma$. To see this note that $Ug\subset gF$ and $F\Gamma$ is closed, so $\cl{Ug\Gamma}\subset gF\Gamma$, and by the orbit closure theorem (cf.~Lesigne~\cite[Theorem 2]{Lesigne}), $\cl{Ug\Gamma}=Lg\Gamma$ and $g\inv L g\in \cH$. Therefore $g\inv Lg\subset F$. 
Since $g\in X(U,g\inv L g)$ and $g\not\in S(U,F)$, we conclude that $\dim (g\inv L g)\geq \dim F$, 
and hence $g\inv L g=F$, and so $\cl{Ug\Gamma}=gF\Gamma$.

Let $\gamma\in \Gamma$. Since 
$g, g\gamma \in G=X(F,U)$, we have that 
$Ug\gamma \subset g\gamma F$, and hence
$gF\Gamma=\cl{Ug\Gamma}\subset g\gamma F\Gamma$.
Therefore 
$gFg\inv \subset g\gamma F \gamma\inv g\inv$, 
and hence $F=\gamma F\gamma\inv$. This proves
that $F$ is normalized by $\Gamma$. 
But $\Ad\Gamma$ is Zariski dense in $\Ad G$ (see~\cite[Chapter 2]{R-book}) and $F$ is
normalized by $G$. 

As we noted, almost every $U$-ergodic component of $\nu$ is of the form $g\mu_F$ for some $g\in G$, and so it is $gFg\inv=F$-invariant. Therefore $\nu$ is $F$-invariant. This completes the proof of the claim.
\end{proof}

We are left with showing that $\nu$ is $G$-invariant.
Let $p\colon G/\Gamma \to G/F\Gamma$ be the natural quotient map, and let $p_{\ast}$ 
denote the corresponding pushforward map of measures.
Then $p_{\ast}(\nu_{t_n}) \to p_{\ast}(\nu)$ as $n\to\infty$. Since $U\subset F$, we have that $\dim((F\Gamma)^0)>\dim(\Gamma^0)$. Therefore by the induction hypothesis stated in the beginning of the proof of Proposition~\ref{prop:local:equi}, $p_{\ast}(\nu)$ is the $G$-invariant probability measure on $G/F\Gamma$. Since $\nu$ is $F$-invariant, it follows that $\nu$ is $G$-invariant.  Since the Haar measure is the unique $G$-invariant probability measure on $X$, we have that  $\nu=\mu$, completing the proof of Proposition~\ref{prop:local:equi}.
\end{proof}

\section{Points with property \texorpdfstring{$\cW(\Gamma^0)$}{W(Γ)}} 
We determine which points $s_0\in (0,1)$ satisfy the property $\cW(\Gamma^0)$. For this purpose, we study the following function using Taylor's expansion and Baker-Campbell-Hausdorff formula: for $u\in (0,1)$ and $\xi\in(-1,1)$ such that $u+\xi\in (0,1)$, we have 
\begin{equation} \label{eq:psi1}
\psi_t(u,\xi)
 =\log[\exp(\rho_t\circ\phi(u+\xi))\cdot \exp(-\rho_t\circ\phi(u))] \in \gbar. 
\end{equation}

Let $u\in (0,1)$ be such that $\phi^{(D)}(u)$ exists. By Taylor's formula, for any $\xi\in\R$ such that $u+\xi\in (0,1)$, we have
\begin{equation} \label{eq:Taylor-error}
\phi(u+\xi)=\sum_{l=0}^{D} \frac{1}{l!} \phi^{(l)}(u)\xi^{l}+\varepsilon(u,\xi)\xi^D,
\end{equation}
where $\varepsilon(u,\xi)\to 0$  as $\xi\to 0$.

Since $\gbar^{(k)}\subset \gbar^{(k-1)}$ for all  $1\leq k\leq \kappa$, we can choose a subspace $V_k\subset \gbar^{(k-1)}$ such that $\gbar^{(k-1)}=V_k\oplus \gbar^{(k)}$. 
Then 
\begin{equation*}
    \gbar=V_1\oplus \cdots \oplus V_\kappa \text{ and }
    \gbar^{(k)}=V_{k+1}\oplus \cdots \oplus V_\kappa.
\end{equation*}
Let $P_k\colon \gbar\to V_k$ denote the corresponding projection.

By the choice of the integers $D_k$ (the polynomial $Q_k\circ\rho_t$ is of degree at most $D_k$), for all $1\leq k\leq \kappa$ and for each $1\leq i\leq D_k$, there exists a linear map $A_{i,k}\colon \lieg\to V_k$, such that
\begin{equation} \label{eq:Pk-rhot}
P_k\circ\rho_t= \sum_{i=0}^{D_k} t^i A_{i,k}.
\end{equation}

Then
\begin{multline} \label{eq:rhot-phi}
    \rho_t\circ\phi(u+\xi)=\sum_{k=1}^{\kappa} (P_k\circ\rho_t)\circ\phi(u+\xi) \\
    =\sum_{k=1}^\kappa\sum_{i=0}^{D_k} \left(\sum_{l=0}^{D}  t^i\xi^l \frac{1}{l!}  A_{i,k} \phi^{(l)}(u) + t^i\xi^{D} A_{i,k}\varepsilon(u,\xi)\right). 
\end{multline}
By the Baker-Campbell-Hausdorff formula, since $\gbar^{(\kappa)}=0$, using the definition of $\psi_t$ in~\eqref{eq:psi1}, we have 
\begin{equation}  \label{eq:psiZ}
\psi_{t}(u,\xi)
=\sum_{i=0}^{\kappa D} \sum_{l=0}^{\kappa D} t^{i}\xi^{l} Z_{(i,l)}(u,\xi),
\end{equation}
where the functions $Z_{(i,l)}$ can be expressed as \begin{equation} \label{eq:ZY}
Z_{(i,l)}(u,\xi)=Y_{(i,l)}(u)+\varepsilon_{(i,l)}(u,\xi),   
\end{equation} 
where $Y_{(i,l)}(u)$ is a fixed linear combination of nested commutators of the form
\begin{equation}
\label{eq:comm}
[X_1,[X_2,[\cdots,[X_{n-1},X_n]]]], \quad\text{ where $1\leq n\leq \kappa$.}
\end{equation}
and
\begin{equation} \label{eq:Xm}
X_m=A_{i_{m},k_{m}}\phi^{(l_m)}(u)\in V_{k_m}, \quad\text{ for $1\leq m\leq n$},
\end{equation}
and $\varepsilon_{(i,l)}(u,\xi)$ is a fixed linear combination
of similar nested commutators where one or more of the $X_m$'s with $l_m=D$ are replaced by
\begin{equation} \label{eq:tildeXm}
\tilde X_m=A_{i_{m},k_{m}}\varepsilon(u,\xi)\in V_{k_m}.
\end{equation}

\begin{Claim} \label{claim:psi}
We have the following:
\begin{enumerate} 
    \item \label{itm:epsilon-ri}
    $\lim_{\xi\to 0} \varepsilon_{(i,l)}(u,\xi)=0$. 
    \item \label{itm:Z0} 
    For all $i$, $Z_{(i,0)}=0$. 
    \item \label{itm:lD}
    If $l<D$, then $\varepsilon_{(i,l)}=0$.
    \item \label{itm:Zil0}
    If $i>D$, then $Z_{(i,l)}=0$.
    \end{enumerate}
In  particular, 
\begin{align} 
\psi_{t}(u,\xi) 
=\sum_{(i,l)\in\cP} t^{i}\xi^{l} Y_{(i,l)}(u) +  \sum_{i=0}^{D} \sum_{l=D}^{\kappa D} t^i\xi^l \varepsilon_{(i,l)}(u,\xi),
\label{eq:Psi_T}
\end{align}
where
\begin{equation} \label{eq:cP}
    \cP=\{(i,l)\colon 0\leq i\leq D,\, 1\leq l\leq \kappa D\}.
\end{equation}
\end{Claim}

\begin{proof}
From the description of $\varepsilon_{(i,l)}(u,\xi)$ and~\eqref{eq:tildeXm}, due to~\eqref{eq:Taylor-error} we conclude that~\eqref{itm:epsilon-ri} holds. 
Since $\psi_t(u,0)=0$, we have that~\eqref{itm:Z0} holds. 
   
Using~\eqref{eq:rhot-phi}, the following relations hold between various indices  in~\eqref{eq:Xm} and~\eqref{eq:tildeXm}:
\begin{equation} \label{eq:dm}
    1\leq k_m\leq \kappa, \ 
    0\leq i_m\leq D_{k_m}, \text{ and } 
    0\leq l_m\leq D.
\end{equation}
In view of equations~\eqref{eq:psiZ}-\eqref{eq:tildeXm}, we obtain
\begin{equation} \label{eq:il}
    \sum_{i=1}^n i_m=i \text{ and }\sum_{m=1}^n l_m=l.
\end{equation}

For every commutator appearing in the expression for $\varepsilon_{(i,l)}(u,\xi)$, we have $l_m=D$ for at least one $m$, and hence $l\geq D$ by~\eqref{eq:il}. Therefore if $l<D$ then $\varepsilon_{(i,l)}=0$, which proves~\eqref{itm:lD}.

By~\eqref{eq:Xm} $X_m\in V_{k_m}$,  and by~\eqref{eq:tildeXm} $\tilde X_m\in V_{k_m}$. Also
\[
V_{k_m}\subset \gbar^{(k_m-1)} \text{ and } [\gbar^{(k-1)},\gbar^{(k'-1)}]\subset \gbar^{(k+k')-1}\quad\text{ for all } k,k'.
\]
Therefore the nested commutator as in~\eqref{eq:comm} and its analogue involving $\tilde X_m$'s belong to $\gbar^{((\sum_{m=1}^n k_m) -1)}$. Since $\gbar^{(\kappa)}=0$, if the nested commutator is nonzero then 
\begin{equation} \label{eq:km}
\sum_{m=1}^n k_m\leq \kappa.
\end{equation}

By~\eqref{eq:dm} and~\eqref{eq:il}, we have that if the nested commutator as in~\eqref{eq:comm} or its analogue involving $\tilde X_m$'s is nonzero then
\begin{equation} \label{eq:imDm}
i=\sum_{m=1}^n i_m\leq\sum_{m=1}^n D_{k_m}.
\end{equation}
Now in view of~\eqref{eq:comm},  \eqref{eq:km} and~\eqref{eq:imDm}, we recall~\eqref{eq:d}:
\[
D=\max\left\{\sum_{m=1}^n D_{k_m}\colon  \sum_{m=1}^n k_m\leq \kappa,\, 1\leq k_m\leq \kappa, \, 1\leq n\leq \kappa\right\}.
\]
Therefore by~\eqref{eq:imDm}, if $i>D$ then $Z_{(i,l)}=0$. This proves~\eqref{itm:Zil0}
\end{proof}

\begin{Claim} \label{claim:Y_11}
For $0\leq i\leq d_1$, 
\[
dq(Y_{(i,1)}(u))=A_i\phi^{(1)}(u).
\]
\end{Claim}

\begin{proof} First note that $dq=dq\circ P_1$. 
Now
\begin{equation}
    P_1\circ Y_{(i,1)}(u) = A_{i,1}\phi^{(1)}(u), \label{eq:P1Yd1}
\end{equation}
because $l=1$ and by~\eqref{eq:comm} and~\eqref{eq:Xm} we have $n=1$, and hence by~\eqref{eq:comm} $l_1=l=1$.

So by~\eqref{eq:Pk-rhot},
\begin{align}
    dq\circ \rho_t 
    = dq\circ (P_1\circ \rho_t) 
    =dq\circ\Bigl(\sum_{i=0}^{D_1} t^i A_{i,1}\Bigr) =\sum_{i=0}^{D_1} t^i dq\circ A_{i,1}. \label{eq:dqrhotA} 
\end{align}
Therefore by~\eqref{eq:torus}, $d_1\leq D_1$ and
\begin{equation} \label{eq:dqA}
    dq\circ A_{i,1} = 
    \begin{cases}
    A_i & \text{if } 0\leq i\leq d_1\\
    0 & \text{if } d_1<i \leq D_1.
    \end{cases} 
\end{equation}
Hence the claim follows from~\eqref{eq:P1Yd1}. 
\end{proof}

\begin{proposition} \label{prop:J1} 
Given $\varepsilon_{0}>0$, there exists a Borel measurable set $J_{1}\subset (0,1)$ such that $\lambda(J_1)\geq 1-\varepsilon_0$ and the following conditions hold:
\begin{enumerate}
\item \label{itm:unif} 
The derivative $\phi^{(\kappa)}(u)$ exists and is uniformly continuous for $u\in J_1$ and it is bounded on $J_1$. In particular,
for each $(i,l)\in\cP$,  $Y_{(i,l)}(u)$ is uniformly continuous and bounded for $u\in J_1$. 
\item \label{itm:tgt} For every $u\in J_{1}$ and every nontrivial unitary character $\chi$ on the torus $G/[G,G]\Gamma=\bar X$, we have
\[
A_i\phi^{(1)}(u)\not\in\ker d\chi \quad\text{ for some } 0\leq i\leq d_1.
\]
\item  \label{itm:ou} For each $(i,l)\in\cP$, $\varepsilon_{(i,l)}(u,\xi)\to 0$ as $\xi\to 0$, and this convergence is uniform
for $u\in J_{1}$. 
\end{enumerate}

\end{proposition}

\begin{proof}
By our assumption, $\phi^{(\kappa)}(u)$ exists for almost all $u\in (0,1)$, and so  by Lusin's Theorem there exists a compact set $J_2\subset (0,1)$ such that $\phi^{(\kappa)}$ is uniformly continuous on $J_2$ and $\lambda(J_2)\geq 1-\varepsilon_0/2$.

By the condition given in~\eqref{eq:tangent}, since $X^\ast\setminus\{1\}$ is countable, there exists a Borel set $J_3\subset J_2$ such that $\lambda(J_2\setminus J_3)=0$ and for every $\chi\in\bar X^\ast\setminus\{1\}$ and every $u\in J_3$, there exists $1\leq i\leq d_1$ such that $A_i\phi^{(1)}(u)\not\in \ker d\chi$. 

Given $p,q\geq 1$, define the set $J_{p,q}$ to be  
\begin{equation*}
\{u\in J_3\colon \text{ for all } (i,l)\in\cP,\,\abs{\varepsilon_{(i,l)}(u,\xi)}\leq 1/p \text{ for all }\abs{\xi} \leq 1/q\}.
\end{equation*}
Then by part~\eqref{itm:epsilon-ri} of Claim~\ref{claim:psi}, the sets $J_{p,q}$ form a nested sequence of sets growing to $J_3$ as $q\to\infty$. 
Choose, $q_p\geq 1$ such that 
\[
\lambda(J_3\setminus J_{p,q_p})<2^{-p}(\varepsilon_0/2)
\]
and set $J_1=\bigcap_{p=1}^\infty J_{p,q_p}$. Then $\lambda(J_3\setminus J_1)\leq \varepsilon_0/2$, and for every $(i,l)\in\cP$,  $\varepsilon_{i,l}(u,\xi)\to 0$ as $\xi\to 0$ uniformly for $u\in J_1$.  
Also $\lambda(J_1)\geq 1-\varepsilon_0$. 
\end{proof}

\begin{corollary} \label{cor:Lambda:J1}
If $J_1$ satisfies conditions~\eqref{itm:unif},\eqref{itm:tgt} and~\eqref{itm:ou} of Proposition~\ref{prop:J1} for $\Gamma$, then $J_1$ also satisfies the same three conditions when we replace $\Gamma$ by any closed subgroup $\Lambda$ of $G$ containing $\Gamma$ and $\rho_t$ by $dp\circ \rho_t$, where $dp\colon \lieg/\Lie(\Gamma^0)\to \lieg/\Lie(\Lambda^0)$ it the natural quotient map. 
\end{corollary}

\begin{proof}
When $\Lambda$ replaces $\Gamma$, each $Y_{(i,l)}$ gets replaced by $dp(Y_{(i,l)})$ and $\varepsilon_{(i,l)}$ gets replaced by $dp(\varepsilon_{(i,l)})$. Therefore it directly follows that Conditions~\eqref{itm:unif} and~\eqref{itm:ou} hold.

Now $\bar X_1=G/[G,G]\Lambda$ is the abelianization of $G/\Lambda$ with $\liet_1=\Lie(G/[G,G]\Lambda^0)$ being its Lie algebra. Let $dr_1\colon \liet\to\liet_1$ the natural quotient map. Then 
$dq$ gets replaced by $dq_1\colon \Lie(G/\Lambda^0)\to \liet_1$ such that $dq_1\circ dp=dr_1\circ dq$.
Then in view of~\eqref{eq:torus},
\[
dq_1\circ(dp\circ \rho_t)=dr_1\circ (dq\circ \rho_t)=\sum_{i=0}^{d_1} t^i(dr_1\circ A_i).
\]

Thus in case of $\Lambda$ replacing $\Gamma$, we have
that $A_i$ gets replaced by $dr_1\circ A_i$ in Condition~\eqref{itm:tgt}. To verify this condition,   
let $\chi_1\in \bar X_1^\ast\setminus\{1\}$ be given. Put $\chi=\chi_1\circ r_1\in \bar X^\ast\setminus\{1\}$. Then 
    \[
    d\chi_1\circ (dr_1\circ A_i)=d\chi\circ A_i.
    \]
    For any $u\in J_1$, by Condition~\eqref{itm:tgt} for $\Gamma$, pick $0\leq i\leq d_1$ such that $d\chi(A_i\phi^{(1)}(u))\neq 0$. Then $d\chi_1((dr_1\circ A_i)\phi^{(1)}(u))\neq 0$.  
\end{proof}

Let ${\cF}_\Gamma$ denote the collection of all connected normal subgroups $F$ 
of $G$ such that $F\supset \Gamma^0$ and $F\Gamma$ is closed. 
In particular, $F/F\cap\Gamma\cong F\Gamma/\Gamma$ is compact.
By~\cite[Chapter~II]{R-book}, there exists a $\Q$-structure on $G$ such that 
$\Gamma^0$ is a $\Q$-subgroup of $G$ and the image of $\Gamma$ on $G/\Gamma^0$ 
consists of integral points with respect to the $\Q$-structure on the quotient 
algebraic group $G/\Gamma^0$. Moreover for any $F\in \cF_\Gamma$, 
we have that $F/\Gamma^0$ must be an algebraic $\Q$-subgroup of $G/\Gamma^0$. Therefore $\cF_\Gamma$ is countable.  

For any $F\in \cF_\Gamma$, let $\pi_F\colon G\to G/F$ be the natural quotient map and let $d\pi_F\colon \lieg \to \lieg/\lief$ denote its differential, where $\lief$ denotes the Lie algebra of $F$. 

For all $(i,l)\in \cP$ and $F\in\cF_\Gamma$, define
\begin{equation} \label{eq:J_0}
    {K}_{(i,l),F}=\{u\in (0,1) \colon \phi^{(D)}(u) \text{ exists and  } d\pi_F(Y_{(i,l)}(u))=0\}.
\end{equation}

Let $J_1$ be a Borel measurable set which satisfies the conditions~\eqref{itm:unif},\eqref{itm:tgt} and~\eqref{itm:ou} of Proposition~\ref{prop:J1}. 

For any $(i,l)\in\cP$ and $F\in\cF_\Gamma$, let
\begin{multline} \label{eq:SdlF}
S_{(i,l),F}=\{s\in J_1\cap {K}_{(i,l),F}\colon \\
\text{$s$ is not a Lebesgue Density point of $J_1\cap {K}_{(i,l),F}$}\}.
\end{multline}
Then by the Lebesgue Density theorem, $\lambda(S_{(i,l),F})=0$. Let
\begin{equation} \label{eq:S1}
S_\Gamma=\bigcup_{(i,l)\in\cP,\, F\in \cF_\Gamma} S_{(i,l),F}.
\end{equation}
Since $\cP$ is finite and $\cF_\Gamma$ is countable, 
\begin{equation} \label{eq:SJ1}
\lambda(S_\Gamma)=0.
\end{equation}

If $\Lambda$ is a closed subgroup of $G$ containing $\Gamma$, then $\Lambda^0\in\cF_{\Gamma}$. Hence $\cF_{\Lambda}\subset \cF_{\Gamma}$. Therefore using Corollary~\ref{cor:Lambda:J1} and~\eqref{eq:S1}, we have that \begin{equation}
\label{eq:Lambda:S}
S_\Lambda\subset S_\Gamma.
\end{equation}

\begin{proposition} \label{prop:cW}
Let $J_1$ be a Borel set which satisfies all the three conditions of Proposition~\ref{prop:J1}. Let $s_0\in J_1\setminus S_{\Gamma}$. Then $s_0$ has property $\cW(\Gamma^0)$. 
\end{proposition}

\begin{proof}. 
Let
\begin{equation}
\cP^{\prime}=\{(i,l)\in \cP\colon  Y_{(i,l)}(s_{0})\neq 0\}.  \label{eq:Pprime}
\end{equation}
In view of~\eqref{eq:J_0},
\begin{equation} \label{eq:s0K}
s_0\in J_1\cap K, \text{ where } {K}:=\bigcap_{(i,l)\in \cP\setminus\cP'} {K}_{(i,l),\Gamma^0}.
\end{equation}

Set 
\begin{equation} \label{eq:alpha}
\alpha=\max\{i/l\colon (i,l)\in \cP^{\prime}\}.
\end{equation}
We first show that $\alpha\geq 1$. It suffices to show that
\begin{equation} \label{eq:dpiY}
    Y_{(i,1)}(s_0)\neq 0 \quad\text{ for some } 1\leq i\leq d_1.
\end{equation}

Let $\chi$ be a nontrivial character on $G/[G,G]\Gamma=\bar X$. Since $s_0\in J_1$, by condition~\eqref{itm:tgt} of Proposition~\ref{prop:J1} pick $1\leq i\leq d_1$ such that $d\chi(A_i\phi^{(1)}(s_0))\neq 0$. Therefore $A_i\phi^{(1)}(s_0)\neq 0$. Therefore by Claim~\ref{claim:Y_11},
\[
dq(Y_{(i,1)}(s_0))=A_i\phi^{(1)}(s_0)\neq 0.
\]
Therefore by~\eqref{eq:Pprime} and~\eqref{eq:dpiY}, there exists $1\leq i\leq d_1$ such that $(i,1)\in \cP^\prime$. Hence by~\eqref{eq:alpha}
\begin{equation}
\alpha\geq i/1\geq 1.
\end{equation}

Fix $\varepsilon>0$. For every $(i,l)\in \cP\setminus \cP'$, $s_0\in J_1\cap {K}_{(i,l),\Gamma^0}$ by~\eqref{eq:s0K}, and since $s_0\not\in S_\Gamma$, $s_0$ is a Lebesgue density point of $J_1\cap {K}_{(i,l),\Lambda^0}$ by~\eqref{eq:SdlF}. 
For every $n\geq 1$, we can choose $k_n\geq n$ and a compact set $I_{n}\subset [0,1]$ such that if we put $t=t_{k_n}$, then
\[
    s_{0}+s\ell_{t} t^{-1}\in J_1\cap {K}_{(i,l),\Gamma^0}\text{ for all $s\in I_{n}$ and $(i,l)\in \cP\setminus \cP'$},
\]
and 
\[
\lambda([0,1]\setminus I_{n})\leq 2^{-n}\varepsilon.
\]
Let $I_{\varepsilon}=\bigcap_{i=1}^{\infty} I_{n}$. 
Then $\lambda((0,1) \setminus I_{\varepsilon})\leq \varepsilon$. 
For any $s\in I_{\varepsilon}$, $\zeta\in\R$, and $n\geq 1$, put
\begin{equation} \label{eq:tuxi}
t=t_{k_n}, \quad u=s_{0}+s\ell_{t}t^{-1}\in J_1\cap{K}, \text{  and } \xi=\zeta t^{-\alpha}.
\end{equation}
Then the following statements hold:
\begin{enumerate}
\item For $l>i$, since $\alpha\geq 1$, by~\eqref{itm:unif} of Proposition~\ref{prop:J1},
\begin{equation}
    t^i\xi^l Y_{(i,l)}(u)\to 0\quad\text{ as $n\to\infty$.}
\end{equation}
\item For $l\geq D\geq i$, since $\alpha\geq 1$ and $u\in J_1$,  by~\eqref{itm:ou} of Proposition~\ref{prop:J1}, 
\begin{equation}
    t^{i}\xi^{l}\varepsilon_{(i,l)}(u,\xi)=\zeta^{l}t^{-(\alpha l-i)}\varepsilon_{(i,l)}(u,\xi) \to 0 
    \quad\text{ as $n\to\infty$.}
\end{equation}
\item For all $(i,l)\in \cP\setminus \cP'$, $u\in K\subset {K}_{(i,l),\Gamma^0}$, so by~\eqref{eq:J_0}, 
\begin{equation}
Y_{(i,l)}(u)=0.
\end{equation}
\item Since $s_{0}, u\in J_1$, for each $(i,l)\in \cP'$ we have that 
\begin{equation}
    Y_{(i,l)}(u)\to Y_{(i,l)}(s_{0}) \quad\text{ as $n\to\infty$,}
\end{equation}
by continuity of $Y_{(i,l)}$ on $J_{1}$.
\item For each $(i,l)\in \cP'$,
\begin{equation}
\lim_{n\to\infty} t^i\xi^l= \lim_{n\to\infty} \zeta^l t^{i-\alpha l}=
\begin{cases} 
\zeta^l, & \text{if $i/l=\alpha$}\\
0, & \text{if $i/l<\alpha$.} 
\end{cases}
\end{equation}
\end{enumerate}

In view of~\eqref{eq:Psi_T}, by the above list of observations
\begin{align} \label{eq:main}
\lim_{n\to\infty} \psi_{t}(u,\xi)
& = \sum_{l\in\cL} \zeta^{l} Y_{(\alpha l,l)}(s_{0}) :=\eta(\zeta),
\end{align}
where $\cL=\{l\colon 1\leq l\leq \kappa D, \, (\alpha l, l)\in \cP'\}$, which is nonempty by the definition of $\alpha$. By~\eqref{eq:Pprime}, 
$Y_{(\alpha l,l)}\neq 0$ for all $l\in\cL$. Therefore $s_0$ has property $\cW(\Gamma^0)$.
 \end{proof}

\begin{proof}[Proof of Theorem~\ref{thm:local:illumination}]
Let $\varepsilon>0$. Obtain $J_1$ as in Proposition~\ref{prop:J1} such that $\lambda(J_1)\geq 1-\varepsilon$. Let $S_\Gamma$ be as defined in~\eqref{eq:S1}. Then by~\eqref{eq:SJ1} we have $\lambda(S_\Gamma)=0$. Let $W_\varepsilon=J_1\setminus S_\Gamma$. Then $\lambda(W_\varepsilon)\geq 1-\varepsilon$. Let $s_0\in W_\varepsilon$. By Corollary~\ref{cor:Lambda:J1} and~\eqref{eq:Lambda:S}, and Proposition~\ref{prop:cW}, $s_0$ has property $\cW(\Lambda^0)$ for every closed subgroup $\Lambda$ of $G$ containing $\Gamma$. Therefore~\eqref{eq:nut} holds by  Proposition~\ref{prop:local:equi}. Now let $W=\cup_{n\in\N} W_{1/n}$. Then $\lambda(W)=1$ and~\eqref{eq:nut} is satisfied for every $s_0\in W$. 
\end{proof}

\section{Weak equidistribution does not imply equidistribution}

\label{sec:counter}

In this section we provide different instances where weak equidistribution of sequence of measures hold, 
but (strong) equidistribution does not. As before let $\lambda$ denote the Lebesgue measure on $\R$ restricted to $(0,1)$.

\begin{proposition}
\label{prop:existweak-not-strong}
 Let $X=\mathbb{T}=G/\Gamma$, where $G=\R$ and $\Gamma=\Z$. Let $\rho_t(v)=tv$ for all $v\in\R=\Lie(G)$ and $t\in\R$, and let $x_{0}\in X$. There exists a measure $\nu$ on $\R$ such that the family of  measures
	$(\mu_t:=\mu_{\nu,x_0,\rho_t})_{t\geq 0}$ is weakly equidistributed but not equidistributed on $(X,\mu)$, where $\mu$ is the Haar measure on $X$.
\end{proposition}
\begin{proof}
We first construct a non-atomic probability measure $\nu$  on $\R$ such $\mu_{3^m}=\mu_{1}$ for all $m\in\N$.  To construct such a measure, let $\psi\colon \R\to \R$ be the function defined by 
	\begin{equation} \label{eq:g}
	\psi(u)=\sum_{\{n\geq 1\colon a_{n}(u)=1\}}\frac{1}{3^{n}}
	\end{equation}
 for 
 $$ u=\sum_{n=1}^{\infty}\frac{a_{n}(u)}{3^{n}} \mod 1,\quad	\text{ where } a_{n}(u)\in\{0,1,2\},  
 $$
and then take  $\nu$ to be the pushforward of $\lambda$ under the map $\psi$ to $\R=\Lie(G)$ and let $\mu_t:=\mu_{\nu,x_0,\rho_t}$. We now check that $\nu$ and $\mu_{t}$ have the stated properties.  

Note that if $0<u_1<u_2<1$ and $\psi(u_1)=\psi(u_2)$ then 
$x_2-x_1=\sum_{n=1}^\infty a_n/3^n$, where $a_n\in\{0,2\}$. So $u_2-u_1$ belongs to the standard Cantor middle third set, which has zero Lebesgue measure. Therefore $\lambda(\psi\inv(\{y\}))=0$ for all $y\in \R$. 
Hence $\nu$ is non-atomic. 

Let $f\in C({\mathbb T})$. Then for all $t\in\R$, 
\begin{equation*}
    \int_X f\,d\mu_t=\int_{0}^{1} f(t\psi(u)+x_0)\,du.
\end{equation*}

For any $u\in [0,1/3)$ and $b\in \{0,1,2\}$,
\begin{equation} \label{eq:3psi}
    3\psi(b/3+u)\equiv \sum_{\{n\geq 1\colon a_n(u)=1\}} \frac{1}{3^{n-1}}
    \equiv \psi(3u)\mod \Z.
\end{equation}
Thus it follows that for any $x_0\in\T$,

\begin{align*}
    \int_0^1 f(3\psi(u)+x_0)\,du
    & =\sum_{b=0}^{2}\int_{b/3}^{b/3+1/3} f(3\psi(u)+x_0)\,du\\ &=\sum_{b=0}^2\int_{0}^{1/3} f(\psi(3u)+x_0)\,du \\
    &=\sum_{b=0}^2 (1/3)\int_{0}^{1} f(\psi(u)+x_0)\,du\\
    &=\int_0^1 f(\psi(u)+x_0)\,du.
\end{align*}
From this we conclude that for any $m\in\N$,
\[
\int_0^1 f(3^m\psi(u)+x_0)\,du=\int_0^1 f(3^{m-1}\psi(u)+x_0)\,du.
\]
It follows that $\mu_{3^{m}}=\mu_1$ for all $m\geq 0$.
Thus the measures $\nu$ and $\mu_t$ satisfy the stated conditions.  
	
 Since $\mu_1$ differs from the Haar measure $\mu$, the family of measures $(\mu_{t})_{t\geq 0}$ is not equidistributed with respect to $\mu$. 
	On the other hand, since $\ker(d\chi)$ is countable set for all nontrivial unitary characters $\chi$, we have that $\nu(\ker(d\chi)+v)=0$ for all $v\in\Lie(G)=\R$, since  $\nu$ is non-atomic. By Theorem~\ref{thm:willumination}, $(\mu_{t})_{t\geq 0}$ is weakly equidistributed.
\end{proof}
This example can be generalized to higher dimensional tori:
\begin{corollary}
	For every torus $X=\mathbb{T}^{d}=G/\Gamma$, where $G=\R^d$, $\Gamma=\Z^d$, let $\rho_t(v)=tv$ for all $v\in\Lie(G)=\R^d$ and $t\in\R^n$. There exists a measure $\sigma$ on $\R^d$ such that for any $x_0\in G/\Gamma$, the family of  measures
	$(\mu_t:=\mu_{\sigma,x_{0},\rho_t})_{t\geq 0}$ is weakly equidistributed but not equidistributed on $(X,\mu)$, where $\mu$ is the Haar measure on $X$.
\end{corollary}
\begin{proof}
	Let $\nu$ be the measure of $\R$ as defined in Proposition~\ref{prop:existweak-not-strong}
and define $\sigma=\nu\times\dots\times\nu$, $\mu_t:=\mu_{\sigma,x_{0},\rho_t}$. We claim that $\sigma$ and $\mu_t$ satisfies these conditions. 
	
	Suppose that $x_{0}=(y_{1},\dots,y_{d})$. Then 
	$$
	\mu_{3^{n}}=\mu_{\nu,y_{1},\rho_{3^{n}}}\times\dots\times\mu_{\mu,y_{d},\rho_{3^{n}}}= \mu_{\nu,y_{1},\rho_{3}}\times\dots\times\mu_{\nu,y_{d},\rho_{3}}=\mu_{3},
	$$
and so $(\mu_{t})_{t\geq 0}$ is not equidistributed.
	
On the other hand,  for every nontrivial unitary character $\chi$ on $\mathbb{T}^{d}$ and every $v\in\Lie(G)=\R^d$, there exist finitely many linear functionals $g_{1},\dots,g_{k}\colon \R^{d-1}\to \R$ and some $1\leq i\leq d$ such that
	$$v+\ker(d\chi)=\bigcup_{m=1}^{k}\{(y_{1},\dots,y_{d})\in \
	\R^{d}\colon g_{m}(y_{1},\dots,\hat{y_{i}},\dots,y_{d})=y_{i}\}.$$
	Assuming without loss of generality that $i=d$, for all $1\leq m\leq k$,
	\begin{equation}\nonumber
	  \begin{split}
	    &\quad\sigma(\{(y_{1},\dots,y_{d})\colon g_{m}(y_{1},\dots,y_{d-1})=y_{d}\})
	    \\&=\int_{\R^{d-1}}\nu(\{y_d\in\lieg\colon y_{d}=g_{m}(y_{1},\dots,y_{d-1})\})\, d\nu(y_{1})\dots d\nu(y_{d-1})
	    \\&=0,
	  \end{split}
	\end{equation}
	because $\nu$ is non-atomic.
	Thus $\sigma(v+\ker(d\chi))=0$ and by Theorem~\ref{thm:willumination}, $(\mu_t)_{t\geq 0}$ is weakly equidistributed.
	\end{proof}

\begin{theorem}
	Let $X=\mathbb{T}^{2}=G/\Gamma$, where $G=\mathbb{R}^{2}$ and $\Gamma=\mathbb{Z}^{2}$ and let $\rho_t(v)=tv$ for all $v\in \Lie(G)=\R^2$ and $t\in\R$.  There is
a function $\phi\colon(0,1)\to\lieg$ such that for any $x_0\in\T^2$, the family of measures $(\mu_t:=\mu_{\phi,x_{0},\rho_t})_{t\geq 0}$ is weakly equidistributed but not equidistributed on $X$.
\end{theorem}

\begin{proof}
	Let $\phi\colon (0,1)\to \Lie(G)=\R^2$ be  defined by
$$	    
\phi(u)=(u,\psi(u))\quad \text{ for all } u\in (0,1),
$$
where $\psi\colon (0,1)\to\R$ is defined by~\eqref{eq:g}. Let  $\mu_{t}=\mu_{\phi,x_0,\rho_{t}}$. 
	
We first show that $(\mu_{t})_{t\geq 0}$ is weakly equidistributed. 

Let $\nu$ be the probability measure on $\R^2$ which is the pushforward of the Lebesgue measure on $(0,1)$ under the map $\phi$. Let \[
C=\{(u,\psi(u))\in \R^2\colon u\in (0,1)\}.
\]
Let $P_1\colon \R^2\to\R$ denote the projection on the first factor. Then for any $E\subset \R^2$ we have $\nu(E)=\lambda(P_1(C\cap E))$. Then $\mu_{\phi,x_0,\rho_t}=\mu_{\nu,x_0,\rho_t}$ for all $t\in\R$. Therefore by Theorem~\ref{thm:willumination}, to prove the weak equidistribution of $(\mu_{\nu,x_0,\rho_t})_{t\geq 0}$, it suffices to show that for all $(p,q)\in\Z^2\setminus\{(0,0)\}$ and $z\in\R$, if we put 
\[
L(p,q,z)=\{(x,y)\colon px+qy=z\}
\]
then 
\[
\nu(L(p,q,z))=\lambda(P_1(L(p,q,z)\cap C))=0.
\]

First suppose that $p=0$. Then $q\neq 0$ and $P_1(L(p,q,z)\cap C)=\psi\inv(z/q)$, which is contained in a translate of a Cantor middle set. Therefore $\lambda(P_1(L(p,q,z)\cap C))=0$. 

Now assume that $p\neq 0$. For $N\in\mathbb{N}$ and
	$a_{1},\dots,a_{N}\in\{0,1,2\}$, set
	\begin{equation*}
	B_{a_{1},\dots,a_{N}}
	=\left[
	\sum_{n=1}^{N}\frac{a_{n}}{3^{n}},\sum_{n=1}^{N}\frac{a_{n}}{3^{n}}+\frac{1}{3^{N}}
	\right)
	\times 
	\left[
	\sum_{a_{n}=1}\frac{1}{3^{n}},\sum_{a_{n}=1}\frac{1}{3^{n}}+\frac{1}{2\cdot 3^{N}}
	\right].
	\end{equation*}
	 Note that for $b_n\in\{0,1,2\}$ for all $n\geq N+1$, 
	\[
	 u = \sum_{n=1}^{N}\frac{a_{n}}{3^{n}}+\sum_{n=N+1}^\infty \frac{b_n}{3^{N}}
	\implies 
	\sum_{a_{n}=1}\frac{1}{3^{n}}\leq \psi(u)
	\leq 
	\sum_{a_{n}=1}\frac{1}{3^{n}}+\sum_{n=N+1}^\infty\frac{1}{3^{n}}.
	\]
	Also 
	\[
	(0,1)\subset\bigcup_{a_1,\ldots,a_N\in\{0,1,2\}} \left[
	\sum_{n=1}^{N}\frac{a_{n}}{3^{n}},\sum_{n=1}^{N}\frac{a_{n}}{3^{n}}+\frac{1}{3^{N}}
	\right).
	\]
	
	Therefore 
	$$C\subset B_{N}=\bigcup_{a_{1},\dots,a_{N}\in\{0,1,2\}}B_{a_{1},\dots,a_{N}}.$$
	
	Therefore it suffices to prove that
	\[
	\lim_{N\to\infty} \lambda(P_1(L(p,q,z)\cap B_N))=0.
	\]
	
	It is easy to see that there exist $s_{1},\dots,s_{2^{N}}\in (0,1)$ such that for all $a_{1},\dots,a_{N}\in\{0,1,2\}$, $\sum_{a_{n}=1}\frac{1}{3^{n}}=s_{i}$ for some $1\leq i\leq 2^{N}$. Let
	\[
		C_{j,N}=
	\bigcup_{a_{1},\dots,a_{N}\in\{0,1,2\},\sum_{a_{n}=1}
	\frac{1}{3^{n}}=s_{j}}B_{a_{1},\dots,a_{N}}.
	\]
	Then
	$$
	B_{N}=\bigcup_{j=1}^{2^{N}}C_{j,N}.
	$$

    Now 
	\[
	P_1(L(p,q,z)\cap C_{j,N})\subset \{(z/p-(q/p)y\colon y\in [s_j, s_j+1/(2\cdot{3^N})]\}.
	\]
	Therefore 
	\[
	\lambda(P_1(L(p,q,z)\cap C_{j,N}))\leq \abs{q/p}/(2\cdot 3^N).
	\]
	Hence
\[
\lambda(P_1(L(p,q,z)\cap B_{N}))\leq 2^N\cdot \abs{q/2p}/3^{N}\to 0 \text{ as $N\to\infty$}.
\]
Thus we can now conclude that the family of measures $(\mu_t)_{t\geq 0}$ is weakly equidistributed.  

	We now prove that $(\mu_{t})_{t\geq 0}$ is not equidistributed. By definition, for any continuous function $f$ on $\T^2$, we have
	$$\int_{\T^2} f\,d\mu_{t}=\int_{0}^{1} f((tu,t\psi(u))+x_0)\,du.$$
	
Recall that by~\eqref{eq:3psi},
For any $u\in [0,1/3)$ and $b\in \{0,1,2\}$,
\begin{equation} 
    3\psi(u+b/3)\equiv \sum_{\{n\geq 1\colon a_n(u)=1\}} \frac{1}{3^{n-1}}
    \equiv \psi(3u) \mod \Z.
\end{equation}
	
Therefore, for any $f\in C(\T^2)$, and $b\in\{0,1,2\}$
	$$
	\int_{\frac{b}{3}}^{\frac{1}{3}+\frac{b}{3}} f((3u,3{\psi}(u))+x_0)\,du=\int_{0}^{1/3} f((3u,\psi(3u))+x_0)\,du.
	$$
	So
	\begin{align*}
\int_{0}^{1} f((3u,3{\psi}(u))+x_0)\,du &=  3\int_{0}^{\frac{1}{3}} f((3u,{\psi}(3u))+x_0)\,du
	\\
	& =  \int_{0}^{1} f((u,{\psi}(u))+x_0)\,du.
	\end{align*}
	From this we can deduce that for any $m\geq 1$,
	\[
	\int_{0}^{1} f((3^mu,3^m{\psi}(u))+x_0)\,du
	=\int_{0}^{1} f((3^{m-1}u,3^{m-1}{\psi}(u))+x_0)\,du.
	\]
	Therefore, 
	$\mu_{3^m}=\mu_{3^{m-1}}=\mu_1$. Since $\mu_1$ is not a Haar measure on $\T^2$, $(\mu_{t})_{t\geq 0}$ is not equidistributed.
\end{proof}

\end{document}